\documentclass[final,leqno,onefignum,onetabnum]{siamltex1213}
\usepackage{epsfig,amsmath,amssymb}
\usepackage{graphpap,latexsym,epsf}
\usepackage{graphics}
\usepackage{color}

\definecolor{green}{rgb}{0.1,0.65,0}
\def\trait #1 #2 #3 {\vrule width #1pt height #2pt depth #3pt}
\def\fin{
      \trait .3 5 0
      \trait 5 .3 0
      \kern-5pt
      \trait 5 5 -4.7
      \trait 0.3 5 0
\medskip}
\newtheorem{teor}{Theorem}[section]
\newtheorem{defin}[teor]{Definition}
\newtheorem{lemm}[teor]{Lemma}
\newtheorem{osse}[teor]{Remark}
\newtheorem{prop}[teor]{Proposition}
\newtheorem{defi}[teor]{Definition}
\newtheorem{coro}[teor]{Corollary}
\newtheorem{prob}[teor]{Problem}
\newtheorem{hypo}[teor]{Hypothesis}

\newcommand{\bele}{\begin{lemm}\begin{sl}}
\newcommand{\enle}{\end{sl}\end{lemm}}
\newcommand{\bedef}{\begin{defi}\begin{sl}}
\newcommand{\eddef}{\end{sl}\end{defi}}
\newcommand{\bete}{\begin{teor}\begin{sl}}
\newcommand{\ente}{\end{sl}\end{teor}}
\newcommand{\beos}{\begin{osse}\begin{rm}}
\newcommand{\eddos}{\end{rm}\end{osse}}
\newcommand{\bepr}{\begin{prop}\begin{sl}}
\newcommand{\empr}{\end{sl}\end{prop}}
\newcommand{\bepro}{\begin{prob}\begin{rm}}
\newcommand{\empro}{\end{rm}\end{prob}}
\newcommand{\bede}{\begin{defin}\begin{sl}}
\newcommand{\edde}{\end{sl}\end{defin}}
\newcommand{\beco}{\begin{coro}\begin{sl}}
\newcommand{\enco}{\end{sl}\end{coro}}
\newcommand{\behy}{\begin{hypo}\begin{sl}}
\newcommand{\enhy}{\end{sl}\end{hypo}}

\newcommand{\duav}[1]{\langle{#1}\rangle}

\newcommand{\RR}{\mathbb{R}}

\newcommand{\NN}{\mathbb{N}}

\newcommand{\dx}{\,{\rm d}x}
\newcommand{\dt}{\,{\rm d}t}
\newcommand{\ds}{\,{\rm d}s}
\newcommand{\dy}{\,{\rm d}y}

\def\qed{\ifmmode % if math mode, assume display: omit penalty etc.
  \else \leavevmode\unskip\penalty9999 \hbox{}\nobreak\hfill
  \fi
  \quad\hbox{\hskip.5em\vrule width.4em height.6em depth.05em\hskip.1em}}\def\endproofsym{\qed}
\renewenvironment{proof}[1][Proof]{\trivlist\item[\hskip\labelsep{\hskip0pt
    {\normalfont\scshape#1.}\hskip .321429\parindent}]\ignorespaces}
{\endproofsym\endtrivlist}
\def\endnobox{\def\endproofsym{}\end{proof}\def\endproofsym{\qed}}

\newcommand{\no}{\nonumber}

\newcommand{\duavg}[1]{\left\langle{#1}\right\rangle}

\DeclareMathOperator{\dive}{div}

\let\TeXchi\chi
\def\chi{{\setbox0 \hbox{\mathsurround0pt
$\TeXchi$}\hbox{\raise\dp0 \copy0 }}}

\def\vb{\mathbf v}
\def\vp{\varphi}

\newfont{\ctv}{msam10}

\def\vpl{\vp^{\bf h}}
\def\yl{y^{\bf h}}

\def\fine{\hfill\kern4pt \vrule height4pt depth0pt width4pt }

\numberwithin{equation}{section}

\title{OPTIMAL DISTRIBUTED CONTROL OF A NONLOCAL CONVECTIVE 
CAHN--HILLIARD EQUATION BY THE VELOCITY IN 3D\thanks{The work of E.R. was supported by the FP7-IDEAS-ERC-StG \#256872 (EntroPhase) and by
GNAMPA (Gruppo Nazionale per l'Analisi Matematica, la Probabilit\`a e le loro Applicazioni) of
INdAM (Istituto Nazionale di Alta Matematica).}}

\author{E. Rocca\footnotemark[2]\ \footnotemark[3]
\and J. Sprekels\footnotemark[3]\ \footnotemark[4]}

\begin{document}
\maketitle

\footnotetext[ 2]{Universit\`a di Milano, Dipartimento di Matematica, ViaSaldini 50, 20133 Milano, Italy. (\email{elisabetta.rocca@unimi.it})}
\footnotetext[ 3]{Weierstrass Institute for Applied Analysis and Stochastics,
  Mohrenstrasse 39, 10117 Berlin,
  Germany. (\email{elisabetta.rocca@unimi.it} and \email{juergen.sprekels@wias-berlin.de})}
\footnotetext[ 4]{Humboldt-Universit\"at
 zu Berlin, Institut f\"ur Mathematik, Unter den Linden 6, 10099
 Berlin, Germany.}
%\footnotetext[5]{Support in common for the first and second
%authors.}

\renewcommand{\thefootnote}{\arabic{footnote}}

\slugger{sicon}{xxxx}{xx}{x}{x--x}
\begin{abstract}
In this paper we study a distributed optimal control problem for a nonlocal convective Cahn--Hilliard equation with degenerate mobility and singular potential in three dimensions of space. While the cost functional is of standard
tracking type, the control problem under investigation cannot easily be treated via standard techniques for two
reasons: the state system is a highly nonlinear system of PDEs containing singular and degenerating terms,
and the control variable, which is given by the velocity of the motion occurring in the convective term, is nonlinearly coupled to the state variable. The latter fact makes it necessary to state rather special regularity assumptions for the admissible controls, which, while looking a bit nonstandard, are however quite natural in the corresponding analytical framework. In fact, they are indispensable prerequisites to guarantee the well-posedness of the associated state system. In this contribution, we 
employ recently proved existence, uniqueness and regularity results for the solution to the
associated state system in order to establish the existence of optimal
controls and appropriate first-order necessary optimality conditions
for the optimal control problem. 
\end{abstract}

\begin{keywords}Distributed optimal control, first-order necessary
  optimality 
conditions, nonlocal models, integrodifferential equations,
{convective} 
Cahn--Hilliard equation, phase separation
\end{keywords}

\begin{AMS}49J20, 49J50, 35R09, 45K05, 74N99\end{AMS}

\pagestyle{myheadings}
\thispagestyle{plain}
\markboth{E. ROCCA AND J. SPREKELS}{OPTIMAL CONTROL OF NONLOCAL
  CAHN--HILLIARD EQUATION}

\section{Introduction}\label{intro}
This paper is concerned with the study of a distributed control problem for a 
Cahn--Hilliard type PDE system that may be considered as
a model for an isothermal phase separation of two constituents taking place in 
a fluid flow whose velocity is given. More precisely, we investigate the
case of a {\em nonlocal} Cahn--Hilliard equation with convective term, degenerate mobility and singular potential. In fact, while the standard Cahn--Hilliard equation (cf., e.g.,
 \cite{C,CH,CH2}) is widely used, it seems that a more realistic version of the
Cahn--Hilliard equation can be characterized by a (spatially) nonlocal free energy.
Although the physical relevance of nonlocal interactions
was already pointed out in the pioneering paper \cite{Ro} (see also \cite[4.2]{Em}
and the references therein), the 
isothermal and non-isothermal models containing nonlocal terms have only recently
been studied from the analytical viewpoint (cf., e.g., \cite{BH1, CKRS, gz, GL1, GL2, KRS2} and the references given there). We also remark that
recently  increasing attention has been paid to nonlocal models also from the 
viewpoint of numerics (cf., e.g., \cite{WiseetAl},
\cite{LowWiseetAl}).

The main difference between local and nonlocal models is given by the choice
of the  interaction potential. {Typically, the nonlocal contribution to the free energy has the form} $\,\int_\Omega k(x,y)\,|\varphi(x) - \varphi(y)|^2\dy\,$, with a given
symmetric kernel $k$ defined on  $\Omega\times \Omega$, where $\Omega$ denotes a 
{(sufficiently regular and bounded) domain in $\RR^3$ in which the phase separation
takes place}; its  local Ginzburg--Landau counterpart
is given by $\,(\sigma/2)|\nabla\varphi(x)|^2$,
where the positive parameter $\,\sigma\,$ is a measure for the thickness
of the interface. Here, $\varphi\,$ represents the local concentration of one of the two phases, {which typically attains} values in a bounded interval, say, in $\,[0,1]$. The local potential can be obtained as a formal limit as $\,m\to \infty\,$ from the nonlocal one with the choice $k(x,y)
= m^{5} k(|m(x-y)|^2)$, where $k$ is a nonnegative function with
support in $[0,1]$. This follows from the formula (which was
formally deduced in \cite{KRS})
\begin{align}\nonumber
\int_\Omega m^{5} k(|m(x-y)|^2)\,|\varphi(x) - \varphi(y)|^2\dy
&= \int_{\Omega_m(x)} k(|z|^2) \left|\frac{\varphi\left(x+ \frac{z}{m}\right)
- \varphi(x)}{\frac{1}{m}}\right|^2\,{\rm d}z\\ \nonumber
\stackrel{m\to\infty}{\longrightarrow}&
\int_{\mathbb{R}^3} k(|z|^2)\duavg{\nabla \varphi(x), z}^2\,{\rm d}z
= \frac{\sigma}{2} |\nabla \varphi(x)|^2\,,
\end{align}
for a sufficiently regular $\varphi$, where
$\sigma = 2/3\int_{\mathbb{R}^3} k(|z|^2)|z|^2\,{\rm d}z$ and $\Omega_m(x) =
m(\Omega - x)$. Here we have used {that}
$\int_{\mathbb{R}^3} k(|z|^2)\duavg{e,z}^2\,{\rm d}z = 1/3\,\int_{\mathbb{R}^3}
k(|z|^2)|z|^2\,{\rm d}z$ for every unit vector $e \in\mathbb{R}^3$.
As a consequence, the local Cahn--Hilliard equation  can be viewed as
an approximation of the nonlocal one and vice versa.
We remark at this point that typical integral kernels, which arise
in applications and meet the regularity assumptions stated below in Section 2, are 
given by the classical Newton potential 
$$k(x)=\kappa\,|x|^{-1}, \quad x\neq 0, \quad\mbox{where }\,\kappa>0 \,\,\,
\mbox{is a constant},$$
by the usual mollifiers, and by the Gaussian kernels 
$$k(x)=\kappa_2\,\exp\left(-|x|^2/{\kappa_3}\right), \quad x\in \RR^3, \quad\mbox{where}
\,\,\,\kappa_2>0\,\,\,\mbox{and} \,\,\,\kappa_3\,\,\,\mbox{are constants}.$$ 

In the seminal paper \cite{EG}, the authors established the existence of a weak solution to the local Cahn--Hilliard equation with degenerate mobility and singular potentials endowed with no-flux boundary conditions. However, in the local case no uniqueness proof is known in case of degenerate mobility and singular potential. This is one of the main advantages of considering the nonlocal potential{: for the nonlocal Cahn--Hillard system, indeed, in the case of periodic boundary conditions, an existence and uniqueness result was proved in \cite{GL2}. Later, a} more general case was considered in \cite{gz}. More recently, the convergence to single equilibria was studied in \cite{LP,lond-pet} (cf. also \cite{GG4} for further results), and in \cite{fgr} the existence of a global attractor for a convective nonlocal Cahn--Hilliard equation with degenerate mobility and singular potential was proved in the three-dimensional case. {Moreover, for} the 
two-dimensional case  also the {long-time} dynamics of  its coupling with 
the Navier--Stokes equation (the nonlocal version of the so-called {\em H-model}) was analyzed in \cite{fgr}. For this model uniqueness of weak solutions and existence of the global attractor  in two dimensions has been recently proved in \cite{FGG}.

Concerning the problem of deriving first-order necessary optimality conditions for optimal control problems involving {\em local} Cahn--Hilliard equations, we can quote the following references: in \cite{wn99}, the authors studied the case of a polynomially growing potential $f$ (in \eqref{p2}) with constant mobility $m$ in \eqref{p1}, 
while more recently in \cite{hw}  the case of the
double obstacle potential $f=I_{[0,1]}$ in \eqref{p2} with constant
mobility $m$ in \eqref{p1} was investigated;  first-order necessary
  optimality conditions were obtained by means of a {regularization}
procedure.  Moreover, the convective 1D case has been dealt with in \cite{zl13}, and 
the recent paper \cite{zl14} discusses the 2D case, where the boundary 
conditions $\varphi=\Delta \varphi=0$ were prescribed in place of the usual no-flux
conditions for $\varphi$ and the chemical potential. Notice that in all of the abovementioned
contributions a distributed control was assumed which was not related to the 
fluid velocity.
Let us finally recall the papers \cite{ColliGilSpr} and \cite{ColliGilPodioSpr}, where the authors studied the optimal control problem associated with a non-standard phase field model of Cahn--Hilliard type, and \cite{blanketal}, respectively, where optimization techniques were used in order to solve variational inequalities related to
 Allen--Cahn and Cahn--Hilliard equations. 

While optimal control problems for certain classes of PDEs coupled with nonlocal 
{\em boundary conditions} have already been studied in the literature (cf., e.g., \cite{dksty11,my09,mpt06,ph}),  
to our best knowledge no analytical contribution exists in  the literature to the study of optimal control problems for nonlocal phase field models of convective Cahn-Hilliard type and, more generally, for nonlocal PDEs where the nonlocal operator appears in the PDEs and not on the boundary. 

Another novelty of this paper is the use of the fluid velocity field as the control parameter. This entails that through the convective term there arises a nonlinear coupling between control and state in product form that renders the analysis difficult. Practical applications of this concept arise (at least indirectly) in the growth of bulk semiconductor crystals. A typical case is the block solidification of large silicon crystals for photovoltaic applications: in this industrial process a mixture of several species of atoms (inpurities) dissolved in the  silicon melt has to be moved by the flow 
(i.e., by the velocity field ${\bf v}$) to the boundary of the solidifying silicon in order to maximize the purified high quality part of the resulting silicon ingot. In other words, the flow pattern acts as a control to optimize the final distribution of the impurities. Notice that in this application the control through the velocity ${\bf v}$ is only {\em indirect}, since the flow pattern is itself controlled via magnetic fields that induce a Lorentz force in the electrically conducting silicon melt. For a  description of such a block solidification process we refer to, e.g., \cite{Kudla}.    

Throughout this paper, we will generally assume that $\Omega\subset \RR^3\,$ is a bounded 
and connected domain with smooth boundary $\,\partial\Omega\,$ and outward unit normal ${\bf n}$, and we denote $Q:=\Omega\times(0,T)$ and $\Sigma:=\partial\Omega\times (0,T)$, where $T>0$ is a prescribed final time. 
We then consider the following control problem:

\vspace{3mm}
\noindent
{
{\bf (CP)} \,\,Minimize the cost functional} 
\begin{align}\label{Ji}
&J(\vp, \vb)=\frac{\beta_1}{2}\int_0^T\!\!\int_\Omega |\vp-\vp_Q|^2\dx\dt
+\frac{\beta_2}{2}\int_\Omega |\vp(T)-\vp_\Omega|^2\dx+\frac{\beta_3}{2}\int_0^T\!\!\int_\Omega |\vb|^2
\dx\dt\,,
\end{align}
subject to 
the {initial-boundary value problem (the {\em state system})}
\begin{align}\label{p1}
&\varphi_t-\dive \left(m(\varphi)\nabla\mu\right)=-\vb\cdot\nabla \varphi\quad\hbox{in }\, Q\,,\\[2mm]
\label{p2}
&\mu=f'(\varphi)+w\quad\hbox{in }\, Q\,,\\[2mm]
\label{p3}
&{w(x,t)=\int_\Omega k(|x-y|)(1-2\varphi(y,t))\dy}\quad\hbox{in }\, Q\,,\\[2mm]
\label{p4}
&m(\varphi)\nabla\mu\cdot {\bf n}=0\quad\hbox{on }\, \Sigma\,,\\[2mm]
\label{p5}
&\vp(0)=\vp_0\quad\hbox{in }\Omega\,,
\end{align}
and to the constraint that the velocity $\vb$, which plays the role of the control, 
belongs to a suitable closed, bounded and convex subset (to be specified later) of the space
\begin{equation}
\label{defVau}
{\cal V}:=\{\vb\in L^2(0,T;H_{div}^1(\Omega))\cap L^\infty(Q)^3:\,\,\,\exists\,\vb_t\in L^2(0,T;L^3(\Omega)^3)\},
\end{equation}
where
\begin{equation}
\label{defHdiv}
H^1_{div}(\Omega):=\{\vb\in H_0^1(\Omega)^3: \mbox{div}(\vb)=0\}\,.
\end{equation}
Notice
 that the velocity is assumed divergence free, and we recall that through the 
convective term $\,-\vb\cdot\nabla \vp\,$ the coupling between control and state is nonlinear.
This nonlinear coupling between control and state is the reason for the strong and a bit
nonstandard regularity assumption for the time derivative of the control $\,\vb$. We also
remark that both $\,H^1_{div}(\Omega)\,$ and $\,{\cal V}\,$ are Banach spaces when
equipped with their natural norms, and that the embedding 
${\cal V}\subset C^0([0,T];L^3(\Omega)^3)$ is continuous. 

The singular potential $f$ will be taken in the typical logarithmic  form  {(cf. the original paper \cite{CH})}
 \[
f(\varphi)=\vp\log(\vp)+(1-\vp)\log(1-\vp),\]
and the mobility $m$, which degenerates at the pure phases $\vp=0$ and $\vp=1$,  has to satisfy the compatibility condition (cf. \cite{EG}, \cite{gz}, \cite{lond-pet})
\[
m(\varphi)=\frac{c_0}{f''(\vp)}=c_0\vp(1-\vp)\, {\quad\mbox{with some constant }\,c_0>0,}
\]
which entails that we have the relations
\begin{equation}\label{flux}
m(\vp)f''(\vp)\equiv c_0\,,\quad\,m(\vp)\nabla\mu=c_0\,\nabla\vp+m(\vp)\,\nabla w.
\end{equation}

\noindent
Moreover, throughout this paper we assume that the given constants  $\beta_1, \,\beta_2,\,\beta_3$ in \eqref{Ji} are nonnegative, while $\vp_Q\in L^2(Q)$ and $\vp_\Omega\in L^2(\Omega)$ represent prescribed target functions of the cost functional $J$.
We could generalize both the expressions of $J$ and of the potential $f$, but we restrict ourselves to the above situation for the sake of a simpler exposition. In particular, we could consider {the case when}
\begin{align}\nonumber
&f\in C^4(0,1) \hbox{ is strictly convex in }\, (0,1), \quad Im(f')^{-1}=[0,1], \\[1mm]
\nonumber
&\frac{1}{f''} \hbox{ is strictly concave in }\,(0,1),
\end{align}
and, for example, 
\[
m\in C^{2}([0,1])\quad\hbox{{satisfies }}\,\,m(\vp)f''(\vp)\geq c_0>0\hbox{ \,for every\, }\vp\in [0,1].
\] 
\noindent
Other interesting problems would be related to the case of more general potentials and mobilities, but also to the optimal control problem related to the coupling of 
(\ref{p1})--(\ref{p5}) with a Navier--Stokes system governing the evolution of the velocity 
$\vb$. The existence of weak solutions {to such coupled systems and their} long-time behavior {have} recently been studied in \cite{fgr} in the {two- and three-dimensional} cases. The analysis of {an} associated control problem in the 2D case will be the subject of a forthcoming paper. 

\paragraph{Plan of the paper} The paper is organized as follows: in Section~\ref{wellposed}, we recall known results regarding the well-posedness of the PDE system (\ref{p1})--(\ref{p5}) as well as the related separation property. We also prove a continuous dependence result (Lemma 2.1) which is needed for the analysis of the control problem. In  Section~\ref{optimalcontrol}, we prove the main results of this paper concerning existence and  
first-order necessary optimality conditions for the optimal control problem {\bf (CP)}.

Throughout this paper we will denote the norm of a Banach space $E$ by $\|\,\cdot\,\|_E$.
In the following, we will make repeated use of Young's inequality
\begin{equation}
\label{young}
a\,b\le \delta \,a^2 + \frac 1{4\delta} b^2 \quad\mbox{for all }\,
a,b\in \RR \,\quad\mbox{and }\, \delta>0,
\end{equation}
as well as of the fact that for three dimensions of space the embeddings $\,H^1(\Omega)\subset L^p(\Omega)$, $1\le p\le 6$, and $\,H^2(\Omega)\subset C^0(\overline{\Omega})$ are continuous and 
(in the first case only for $1\le p<6$) 
compact. Moreover, we recall that for smooth and bounded three-dimensional domains there hold the special
Gagliardo--Nirenberg inequalities
\begin{eqnarray}
\label{GN}
\|v\|_{L^3(\Omega)}&\le&\widehat{K}_1\left(\|v\|_{L^2(\Omega)}^{1/2}\,\|v\|_{H^1(\Omega)}^{1/2} \,+\,\|v\|_{L^2(\Omega)}\right)\quad\,\,\forall\,v\in H^1(\Omega)\,,\\[1mm]
\label{GN2}
\|v\|_{L^4(\Omega)}&\le&\widehat{K}_2\left(\|v\|_{L^2(\Omega)}^{1/4}\,\|v\|_{H^1(\Omega)}^{3/4} \,+\,\|v\|_{L^2(\Omega)}\right)\quad\,\,\forall\,v\in H^1(\Omega)\,,
\end{eqnarray} 
where the constants $\widehat{K}_1>0$ and $\widehat{K}_2>0$ depend only on $\Omega$; 
observe that (\ref{young}) and the continuity of the embedding $W^{1,4}(\Omega)\subset L^\infty(\Omega)$ imply that for every $\delta>0$ it holds
\begin{eqnarray}
\label{compactness}
\|v\|_{L^3(\Omega)}^2&\le&\delta\,\|v\|_{H^1(\Omega)}^2\,+\,\frac {\widehat{K}_ 3} \delta 
\,\|v\|_{L^2(\Omega)}^2 \quad\,\,\forall\,v\in H^1(\Omega)\,,\\[1mm]
\label{compactness2}
\|v\|_{L^\infty(\Omega)}^2&\le&\delta\,\|v\|_{H^2(\Omega)}^2\,+\,\frac {\widehat{K}_ 4} \delta 
\,\|v\|_{H^1(\Omega)}^2 \quad\,\,\forall\,v\in H^2(\Omega)\,,
\end{eqnarray}
where also $\widehat{K}_3>0$ and $\widehat{K}_4>0$ depend only on $\Omega$. 
We also recall the 
well-known fact that the trace operator $\varphi\mapsto\varphi_{|\partial\Omega}$ is 
a continuous mapping from $ H^1(\Omega)\cap
L^\infty(\Omega)$ into $H^{1/2}(\partial\Omega)\cap L^\infty(\partial\Omega)$; moreover, it follows from the form of the intrinsic norm of $\,H^{1/2}(\partial
\Omega)\,$ that we have for products the implications
\begin{align}
\label{product1}
u, v\in H^{1/2}(\partial\Omega)\cap L^\infty(\partial\Omega)
&
\quad \Longrightarrow
\quad u\,v \in H^{1/2}(\partial\Omega)\cap L^\infty(\partial\Omega),
\\
\label{product2}
u,v \in L^2(0,T;H^{1/2}(\partial\Omega))\cap L^\infty(\Sigma) 
&
\quad
\Longrightarrow \quad u\,v \in L^2(0,T;H^{1/2}(\partial\Omega)\cap L^\infty
(\partial\Omega))\,.
\end{align}

Finally,
for the sake of a shorter exposition, we 
denote by $\,{\cal K}\,$ the integral operator that assigns to $\vp$ the function $w$ through (\ref{p3}); that is, we put
\begin{equation}
\label{intop}
{\cal K}(\vp)(x,t):=\int_\Omega k(|x-y|)(1-2\vp(y,t))\dy\,.
\end{equation}

\section{Well-posedness {of the state system}}
\label{wellposed}

In the following, we study the state system (\ref{p1})--(\ref{p5}). 
To fix things, we assume for the set of
admissible controls:

\vspace{5mm}\noindent
{\bf (H1)} \quad\,${\cal V}_{\rm ad}\,:=\,\left\{\vb=(v_1,v_2,v_3)
\in {\cal V}: \,\,\,\widetilde v_{1_i}
\le v_i\le \widetilde v_{2_i}\,\,\,\mbox{a.e. in}\,\,\, Q,\,\,\, i=1,2,3,
\right.$\\[1mm]
$\hspace*{27mm}\left. \|\vb\|_{L^2(0,T;H^1(\Omega)^3)}\,+\,\|\vb_t\|
_{L^2(0,T;L^3(\Omega)^3)}\,\le\,V \right\},$\\[2mm]
where $\,V>0\,$ is a given constant and $\,\widetilde v_{1_i}, 
\widetilde v_{2_i}\in L^\infty(Q)$,
$i=1,2,3$,  are given threshold functions; we generally assume
  that ${\cal V}_{\rm ad}\not= \emptyset$.

\noindent
Observe that $\,{\cal V}_{\rm ad}\,$ is a bounded, closed, and convex subset of 
$\,{\cal V}$, which is
certainly contained in some bounded open subset of $\,{\cal V}$. For convenience, we fix such a set once and for all,
noting that any other such set could be used instead:

\vspace{3mm}
\noindent
{\bf (H2)} \,\quad 
${\cal V}_R\subset {\cal V}\,$ is an open set satisfying $\,{\cal V}_{\rm ad}\subset
{\cal V}_R\,$ such that, for all $\,\vb\in {\cal V}_R$,
\begin{equation}
\label{defR}
\|\vb\|_{L^2(0,T;H^1(\Omega)^3)}\,+\,\|\vb\|_{L^\infty(Q)^3}\,+\,\|\vb_t\|
_{L^2(0,T;L^3(\Omega)^3)}\,\le\,R\,. 
\end{equation} 

  Before stating some results on the well-posedness of the 
	state system (\ref{p1})--(\ref{p5}), we now formulate the general assumptions for the problem data. We remark at this place that not all of these assumptions are needed
	to ensure the respective results concerning existence, separation, uniqueness, and regularity;
	however, they are indispensable prerequisites for the continuous dependence result of 
	Lemma 2.2 below, which will be needed for the derivation of necessary optimality
	conditions for the control problem. Since we focus on optimal control here, we
	have decided to impose the corresponding (stronger) conditions from the very beginning 
	in order to avoid any confusion. We make the following assumptions:

\vspace{2mm}
\noindent
{\bf (H3)} \,\quad$\vp_0\in H^2(\Omega)$, there is some $\,\kappa_0>0\,$
 such that  \,\,$0 < \kappa_0\le \vp_0\le 1-\kappa_0< 1\,$ \linebreak
\hspace*{13mm} a.e. in $\,\Omega$,
and it holds a.e. in $\,\Omega\,$ that
\begin{align*}
0& = \Big(c_0\,\nabla\varphi_0\,+\,
m(\varphi_0)\,\nabla\int_\Omega k(|x-y|)(1-2\,\varphi_0(y))\dy\Big)
\cdot {\bf n}\\[1mm]
&=m(\vp_0)\,\nabla \mu(\cdot,0) \cdot {\bf n} .
\end{align*}

\vspace{2mm}
\noindent
{\bf (H4)} \,\quad$f(\vp)=\vp\,\log(\vp)+(1-\vp)\,\log(1-\vp)\,$ \,\,for \,$\,0<\vp<1$, 
\,\,\,\,$f(0)=f(1)=0$, \\[1mm]
\hspace*{13mm}$f(\vp_0)\in L^1(\Omega)$.

\vspace{2mm}
\noindent
{\bf (H5)} \,\quad$m(\vp)=\frac{c_0}{f''(\vp)}$ \,for $\,0< \vp < 1$, 
with some $\,c_0>0$. 	
 
\vspace{2mm}
\noindent
{\bf (H6)} \,\quad$\int_\Omega\int_\Omega k(|x-y|)\dx\dy=:k_0<+\infty, \quad 
\sup _{x\in \Omega}\int_\Omega|k(|x-y|)|\dy=:\bar k<+\infty\,.$

\vspace{2mm}
\noindent
{\bf (H7)} \,\quad$\forall\, p\in [1,+\infty]\,\,\,\exists\, k_p>0\,:\,\left\|
-2\int_\Omega k(|x-y|)\,z(y)\dy\right\|_{W^{1,p}(\Omega)}\leq 
{k_p}\, \|z\|_{L^p(\Omega)}\,$\\[2mm]
\hspace*{13mm} for all \,$z\in W^{1,p}(\Omega)$.

\vspace{2mm}
\noindent
{\bf (H8)} \,\quad For $p\in \{2,3\}\,$ there is some $\,s_p>0\,$ 
such that for all $\,z\in W^{1,p}(\Omega)$ \,it holds\\[2mm]
\hspace*{13mm}
$\left\|-2\int_\Omega k(|x-y|)\,z(y) \dy\right\|_{W^{2,p}(\Omega)}\leq s_p\,\|z\|_{W^{1,p}(\Omega)}\,.$

\vspace{5mm}
We now establish some results for the state system. The following result was essentially shown in \cite[Thm. 2.2]{lond-pet} for the case $\vb=0$:

\vspace{5mm}
\noindent
{\sc \quad Proposition 2.1.} \,\,\,{\em The system} (\ref{p1})--(\ref{p5})
{\em admits under the hypotheses} {\bf (H1)}--{\bf (H8)} {\em for any $\,\vb\in {\cal V}_R\,$ a unique solution triple $\,(\varphi,w,\mu)$\, such that}
\begin{equation}
\label{juerg1}
\varphi\in C^1([0,T];L^2(\Omega))\cap H^1(0,T;H^1(\Omega))\cap L^\infty(0,T;H^2(\Omega))\cap
C^0(\overline{Q}).
\end{equation}
{\em Moreover, there is some $\,\kappa\in (0,1)$, which does not depend on the
choice of $\,\vb\in {\cal V}_R$, such that}
\begin{equation}
\label{juerg2}
0<\kappa\le \varphi\le 1-\kappa<1 \quad\,a.e.\,\,\,in\,\,\,Q \,.
\end{equation}

\noindent
{\sc Proof.}  At first, adapting the proof by Gajewski and Zacharias (see \cite[Thm. 3.5]{gz} and also \cite[Thm. 4 and Prop. 4]{fgr}) to the case $\vb\neq 0$,  one can establish the existence of a unique weak solution $\,(\varphi,w,\mu)\,$ to
(\ref{p1})--(\ref{p5}) such that
\begin{align}
\label{juerg3}
&\quad\vp\in H^1(0,T;H^1(\Omega)^*)\cap L^2(0,T;H^1(\Omega))\cap
C^0([0,T];L^2(\Omega)), 
\\[2mm]
\label{juerg4}
&
\quad\int_\Omega\vp(x,t)\dx=\int_\Omega \vp_0(x)\dx \quad\forall\,t\in [0,T],
\quad 0<\vp<1 \quad\mbox{a.e. in }\,Q\,,
\\[2mm] 
\label{juerg5}
&
\quad w\in L^\infty(0,T;W^{1,\infty}(\Omega))\,,
\end{align}
%\\[2mm]
\begin{align}
\label{juerg6}
&
\quad\int_0^T\!\!\int_\Omega m(\vp)\,|\nabla\mu|^2\dx\dt\,<\,+\infty\,.
\end{align}

Next, it is not difficult to see that the additional convective
term $\,-{\bf v}\cdot \nabla\varphi\,$ on the right-hand side of (\ref{p1}) does 
not create major problems in modifying the proof of \cite[Prop. 3.1]{lond-pet}
to the convective case provided the velocity is (as in our case) bounded; in fact, 
just as there it turns out that the expressions $\,\|\ln(\varphi(t)\|_{L^r(\Omega)}\,$ and 
$\,\|\ln(1-\varphi(t))\|_{L^r(\Omega)}\,$ are bounded by a constant that neither depends 
on $r\in [1,+\infty)$ nor on $t\in [0,T]$, whence it can be concluded
that there is a constant $\,\kappa\in (0,1)$, which is independent of the choice of
$\,\vb\in {\cal V}_R$, such that the weak solution satisfies the separation property
(\ref{juerg2}).

In order to prove the regularity property (\ref{juerg1}), we can follow the lines of the proof of \cite[Thm. 2.2]{lond-pet} in which the asserted regularity was shown for the case without convection. We provide here the details
of the argument, since they differ from those given there. To this end, we will first show that (cf. Eq. (4.1) in \cite{lond-pet})
\begin{equation}
\label{new1}
\varphi_t\in L^\infty(0,T;H^1(\Omega)^*)\cap L^2(Q)\,.
\end{equation}

The derivation of (\ref{new1}) requires the introduction of a functional analytic tool which is standard in the framework of Cahn--Hilliard equations. To this end, we denote by $\langle\,\cdot\,,\,\cdot\,\rangle$ the dual pairing between $H^1(\Omega)^*$ and 
$H^1(\Omega)$, and denoting by $\,|\Omega|\,$ the Lebesgue measure of $\Omega$,
we introduce for functions $\psi\in H^1(\Omega)^*$ and $\vp \in L^1(0,T;H^1(\Omega)^*)$ the generalized mean values
\begin{equation}
\label{new3}
\psi^\Omega:=\frac 1 {|\Omega|}\langle \psi, \mathbf{1}\rangle,\quad\mbox{and }\,
\vp^\Omega(t):=(\vp(t))^\Omega \quad\mbox{for a.e. } t\in (0,T).
\end{equation}
We then introduce the operator ${\cal N}$ as the inverse of the Laplacian with zero Neumann boundary condition as follows: we define
$$
{\rm dom}\,{\cal N}:=\left\{\psi_*\in H^1(\Omega)^*: \psi_*^\Omega=0\right\}, \quad\mbox{and }\,
{\cal N}:{\rm dom}\,{\cal N}\to \left\{\psi\in H^1(\Omega): \psi^\Omega=0\right\}
$$
by setting
\begin{equation*}
{\cal N}\psi_*\in H^1(\Omega), \quad ({\cal N}\psi_*)^\Omega=0, \quad\mbox{and }\,
\int_\Omega \nabla{\cal N}\psi_*\cdot\nabla z\dx =\langle\psi_*,z\rangle\quad\forall\,
z\in H^1(\Omega)\,.
\end{equation*}
In other words, $\psi={\cal N}\psi_*$ is the unique solution to the generalized Neumann problem
$\,-\Delta\psi=\psi_*\,$ in $\Omega$, $\partial \psi/\partial{\bf n}=0$ on $\partial \Omega$,
that has zero mean value. It is a well-known fact that through the formula
\begin{equation}
\label{new4}
\|\psi_*\|_*^2\,:=\,\left\|\nabla{\cal N}\left(\psi_*-\psi_*^\Omega\right)\right\|^2
_{L^2(\Omega)}\,+\,\left|\psi_*^\Omega\right|^2\quad\forall\,\psi_*\in H^1(\Omega)^*
\end{equation}
a norm is defined on $H^1(\Omega)^*$, which is equivalent to the standard norm of 
$H^1(\Omega)^*$ and has the following properties:
\begin{eqnarray}
\label{new5}
&&
\,\,\quad\langle \psi_*,{\cal N}\varphi_*\rangle =\langle \varphi_*, {\cal N}\psi_*\rangle
=\int_\Omega (\nabla{\cal N}\psi_*)\cdot (\nabla{\cal N}\varphi_*)\dx \quad
\forall\,\varphi_*,\psi_*\in {\rm dom}\,{\cal N},
\\[1mm]
\label{new6}
&&
\,\,\quad\langle \psi_*, {\cal N}\psi_*\rangle \,=\,\|\psi_*\|_*^2\,=\int_\Omega
|\nabla{\cal N}\psi_*|^2\dx\quad\forall\,\psi_*\in {\rm dom}\,{\cal N},\\[1mm]
\label{new7}
&&
\,\,\quad 2\,\langle\partial_t\psi_*(t), {\cal N}\psi_*(t)\rangle\,=\,\frac {d}{\dt}
\int_\Omega|\nabla{\cal N}\psi_*(t)|^2\dx \,=\,\frac{d}{\dt}\,\|\psi(t)\|_*^2
\quad\mbox{a.e. in } (0,T),\\[1mm]
&&
\,\,\quad\mbox{for any } \psi_*\in H^1(0,T;H^1(\Omega)^*) \,\mbox{ satisfying }
\,\psi_*^\Omega(t)=0 \,\mbox{ for a.e. }\,t\in (0,T).
\nonumber
\end{eqnarray} 
 
We are now in the position to prove (\ref{new1}). In the remainder of the proof,
we will by $C$ denote generic positive constants that depend only on the data of the system and may change within in formulas and/or even within lines. Moreover, we will  argue formally, noting that all of the following arguments can be made rigorous by using difference quotients with respect to time. 

Now  
recall that $\vb(t)$ is divergence free and vanishes on $\partial\Omega$ for almost all
$t\in (0,T)$, whence it follows that $\,\vb(t)\cdot\nabla\vp(t)\in L^2(\Omega)\,$ has zero mean value. It is thus an easy consequence of (\ref{p1}) and (\ref{p4}) that 
$\,\vp_t(t)\,$ belongs to $\,{\rm dom}\,{\cal N}\,$ for almost every $t\in (0,T)$. We may therefore
(formally) differentiate the variational formulation of the state system 
(\ref{p1})--(\ref{p5}) with respect to $t$ and insert ${\cal N}\vp_t(t)$ as test function. As in the proof of \cite[Thm. 2.2]{lond-pet}, this leads for almost
every $t\in (0,T)$ to an estimate of the form
\begin{equation}
\label{new8}
\|\vp_t(t)\|_*^2 \,+\int_0^t\!\!\int_\Omega \vp_t^2 \dx\ds \,\le \,
\|\vp_t(0)\|_*^2\,+\,C \Big (
\int_0^t \|\vp_t(s)\|_*^2\ds \,+\,I_1(t)\,+\,I_2(t)\Big)\,,
\end{equation}
where the terms $I_1(t)$ and $I_2(t)$ originate from the convective term and will be estimated below.

Notice that  $\,\|\vp_t(0)\|_*\,$ is bounded since this is true
for $\,\|\vp_t(0)\|_{L^2(\Omega)}$; indeed, the assumption $\varphi_0\in
H^2(\Omega)$, in combination with (\ref{flux}) and {\bf (H7)}, yields that $\,{\rm div}(m(\vp_0)\nabla \mu(0))\in L^2(\Omega)$, and
since $\vb\in C^0([0,T];L^3(\Omega)^3)$ and $\nabla\vp_0\in L^6(\Omega)$, it is easily seen that also $\,\vb(0)\cdot\nabla\vp_0\in L^2(\Omega)$.

Next, we have, using the fact that $\vb(t)\in H^1_{div}(\Omega)\cap
L^\infty(\Omega)^3$ for almost every $t\in (0,T)$, and invoking (\ref{new6}),
\begin{align}
\label{new9}
I_1(t)\,
&
= \Big|\int_0^t\!\!\int_\Omega (\vb\cdot\nabla\vp_t)\,{\cal N}
\vp_t\dx\ds\Big |
\,=\,\Big|\int_0^t\!\!\int_\Omega \vp_t\, (\vb\cdot\nabla{\cal N}\vp_t)\dx\ds\Big |
\\
&
\le \int_0^t\|\vb(s)\|_{L^\infty(\Omega)^3}\,\|\vp_t(s)\|_{L^2(\Omega)}\,
\|\nabla{\cal N}\vp_t(s)\|_{L^2(\Omega)}\ds
\nonumber\\
&
\le\frac 1 2\int_0^t\!\!\int_\Omega \vp_t^2\dx\ds\,+\,C\int_0^t
\|\vp_t(s)\|_*^2\ds\,. 
\nonumber
\end{align}

\noindent
 Also,
\begin{align}
\label{new10}
I_2(t)\,
&
=\Big|\int_0^t\!\!\int_\Omega (\vb_t\cdot\nabla\vp)\,{\cal N}\vp_t\dx\ds\Big |
\\
&
\le\int_0^t\|\vb_t(s)\|_{L^3(\Omega)^3}\,\|\nabla\vp(s)\|_{L^2(\Omega)^3}\,
\|{\cal N}\vp_t(s)\|_{L^6(\Omega)}\ds
\nonumber\\
&
\le C \,+\,\int_0^t\|\nabla\vp(s)\|^2_{L^2(\Omega)^3}\,\|\varphi_t(s)\|_*^2\ds\,,
\nonumber
\end{align}
where we have used that
\begin{equation*}
\|{\cal N}\vp_t(s)\|_{L^6(\Omega)}\,\le\,C\,\|{\cal N}\vp_t(s)\|_{H^1(\Omega)}
\le\,C\,\|\nabla{\cal N}\vp_t(s)\|_{L^2(\Omega)^3}
\,\le\,C\,\|\vp_t(s)\|_*\,.
\end{equation*}

\noindent
Combining  (\ref{new8})--(\ref{new10}), and noting that the function
$\,s\mapsto \|\nabla\vp(s)\|^2_{L^2(\Omega)^3}\,$ is known to belong to $L^1(0,T)$,
we can finally verify the claim (\ref{new1}) using Gronwall's lemma.

\vspace{2mm}
Next, we can infer from (\ref{flux}), (\ref{new1}), and from the fact that $\,\,-\vb\cdot\nabla\vp\in L^2(Q)$,  that 
$$
\dive (m(\vp)\,\nabla\mu)=c_0\,\Delta \vp \,+\,m'(\vp)\,\nabla\vp\cdot\nabla w
\,+\,m(\vp)\,\Delta w
$$
belongs to $L^2(Q)$. But then it follows from (\ref{juerg3}), (\ref{juerg5})
and {\bf (H8)} that also
\begin{equation}
\label{DeltaL2}
\Delta\vp\in L^2(Q)\,.
\end{equation}

\noindent
Moreover, we know already from (\ref{juerg3}), (\ref{juerg5}), and
{\bf (H8)}, that $\,\nabla w \in L^2(0,T;H^1(\Omega)^3)\linebreak
\cap L^\infty(Q)^3$, so that
$\,\partial w/\partial {\bf n}= (\nabla w)_{|\partial\Omega}\cdot {\bf n}\,$ belongs to $\,L^2(0,T;	H^{1/2}(\partial\Omega))\cap L^\infty(\Sigma)$. Since 
obviously $\,m(\varphi)_{|\partial\Omega}\,$ belongs to the same space, it follows
from the boundary condition (\ref{p4}) and the product rule (\ref{product2}) 
that the same is true for $\,\partial\vp/\partial {\bf n}$. Hence we can infer
from standard elliptic estimates that
\begin{equation}
\label{L2H2}
\vp\in L^2(0,T;H^2(\Omega))\,.
\end{equation}
It then follows from the continuity of the embedding 
$\,H^1(0,T;L^2(\Omega))\cap L^2(0,T;H^2(\Omega))$ \linebreak
$\subset\, C^0([0,T];
H^1(\Omega))\,$ and from {\bf (H7)} that also
\begin{equation}
\label{CH1}
\vp\in C^0([0,T];H^1(\Omega)),\quad\,w\in C^0([0,T];H^2(\Omega)),
\end{equation}
and analogous reasoning as above shows that we also
have 
\begin{equation}
\label{phirand}
\frac {\partial\vp}{\partial {\bf n}}\,\in L^\infty(0,T;H^{1/2}(\partial\Omega))\,.
\end{equation}

\vspace*{2mm}
In the next step we show that it holds (cf. Eq. (4.3) in \cite{lond-pet})
\begin{equation}
\label{phitreg}
\varphi_t\in C^0([0,T];L^2(\Omega))\cap L^2(0,T;H^1(\Omega))\,.
\end{equation}

To this end, we  differentiate 
the variational formulation of problem (\ref{p1})--(\ref{p5}) with respect to time
again and test by $\vp_t$. As in \cite{lond-pet}, we obtain for every $t\in (0,T]$ an inequality of the form
\begin{equation}
\label{new11}
\|\vp_t(t)\|_{L^2(\Omega)}^2\,+\,\int_0^t\!\!\int_\Omega |\nabla\vp_t|^2\dx\ds\,
\le\,\|\vp_t(0)\|_{L^2(\Omega)}^2\,+\,C\int_0^t\!\!\int_\Omega \vp_t^2\dx\ds\,
+\,C\,I(t),
\end{equation}
where we have  $\,\vp_t(0)\in L^2(\Omega)\,$ and where the expression 
$I(t)$ originating from the convective term has to be estimated. Employing 
(\ref{defR}), (\ref{new1}), and (\ref{L2H2}), and invoking 
H\"older's and Young's inequalities as well as the continuity of the embedding $H^1(\Omega)\subset L^6(\Omega)$,  we have, for any $\delta>0$,
\begin{align}
\label{new12}
I(t)
&
=\Big |\int_0^t\!\!\int_\Omega (\vb\cdot\nabla\vp_t+\vb_t\cdot\nabla\vp)\,\vp_t
\dx\ds\Big |
\\
&
\le \int_0^t\|\vb(s)\|_{L^\infty(\Omega)^3}\,\|\nabla\vp_t(s)\|_{L^2(\Omega)^3}\,
\|\vp_t(s)\|_{L^2(\Omega)}\ds\nonumber
\\
&
\quad +\int_0^t \|\vb_t(s)\|_{L^3(\Omega)^3}\,\|\nabla\vp(s)\|_{L^6(\Omega)^3}\,
\|\vp_t(s)\|_{L^2(\Omega)}\ds\nonumber
\\
&
\le \,\frac{C}{\delta}
\,\Big(1\,+\int_0^t\left(1 \,+\,\|\vb_t(s)\|^2_{L^3(\Omega)^3}\right)\,\|\vp_t(s)\|^2
_{L^2(\Omega)} \ds\Big)\nonumber
\\
&
\quad +\,\delta\int_0^t\!\!\int_\Omega |\nabla\vp_t|^2\dx\ds\,.\nonumber
\end{align}
Observing that the function $\,\,s\mapsto\|\vb_t(s)\|^2_{L^3(\Omega)^3}\,\,$
belongs to $L^1(0,T)$, and adjusting $\delta>0$ appropriately small,
we obtain (\ref{phitreg}) by an application of Gronwall's lemma.  

Next, we observe that we have almost everywhere in $Q$ that
\begin{equation}
\label{linfL2}
c_0\,\Delta\vp=\vp_t\,+\,{\bf v}\cdot\nabla\vp\,-\,m'(\vp)\,\nabla\vp
\cdot\nabla w\,-\,m(\vp)\,\Delta w\,,
\end{equation}
and since all terms on the right-hand side are known to belong to 
$L^\infty(0,T;L^2(\Omega))$, the same holds for $\Delta\vp$. Invoking
(\ref{CH1}) and (\ref{phirand}), we therefore obtain from standard
elliptic estimates that
\begin{equation}
\label{linfH2}
\vp\in L^\infty(0,T;H^2(\Omega))\,.
\end{equation}

Finally, we conclude from the continuity of the embedding
$H^1(0,T;H^1(\Omega))\cap L^\infty(0,T;H^2(\Omega))\subset C^0([0,T]; 
H^s(\Omega))$ for every $s\in [0,2)$ that also $\vp\in C^0(\overline{Q})$,
which concludes the proof of the assertion. \qed

\vspace{3mm}
\noindent
{\bf Remark 1.} \,\,A closer inspection of the above proof reveals that
there is a constant $\widehat C_1>0$, which only depends on the data of the system and on the constant $R$, 
such that we have 
\begin{equation}
\label{phibound}
\|\vp\|_{C^1([0,T];L^2(\Omega))\cap H^1(0,T;H^1(\Omega))\cap L^\infty(0,T;H^2(\Omega))\cap
C^0(\overline{Q})}\,\le\,\widehat C_1\,,
\end{equation} 
whenever $\,\varphi\,$ is the first component of a solution $\,(\vp,w,\mu)\,$ associated 
with some $\,\vb\in {\cal V}_R$. But then it follows from the hypotheses {\bf (H7)} and {\bf (H8)} 
that, in particular,
\begin{equation}
\label{wbound}
\|w_t\|_{L^2(0,T;H^2(\Omega))\cap C^0([0,T];H^1(\Omega))}\,+\,\|w\|_{C^0([0,T];H^2(\Omega))\cap L^\infty(0,T;
W^{1,\infty}(\Omega))}\,\le\,
\widehat{C}_2,
\end{equation}
where also $\,\widehat C_2\,$ only depends on the data and $R$. Moreover, the separation
property (\ref{juerg2}) holds even pointwise for every $(x,t)\in \overline{Q}$, whence it follows that
\begin{equation}
\label{f'bounds}
{\max_{1\le i\le 4} \,\|f^{(i)}(\vp)\|_{C^0(\overline{Q})}\,\le\,\widehat C_3,}
\end{equation} 
where, again, $\,\widehat{C}_3\,$ only depends on the data and $R$. Therefore, we can 
conclude from (\ref{p2}) that 
\begin{equation}
\label{mubound}
\|\mu_t\|_{C^0([0,T];L^2(\Omega))\cap L^2(0,T;H^1(\Omega))}\,+\,\|\mu\|_{L^\infty(0,T;H^2(\Omega))}\,\le\,
\widehat{C}_4,
\end{equation}
where also $\,\widehat C_4\,$ only depends on the data and $R$. In the remainder of this paper,
we denote  $\,\,K_1^*\,:=\,\max_{1\le i\le 4} \widehat C_i$.
 
\vspace{3mm}
\noindent
{{\bf Remark 2.} \,\,The separation property (\ref{juerg2}) and the hypotheses {\bf (H4)}
and {\bf (H5)} also entail the estimate}
\begin{equation} \label{j1}
 \frac {c_0} 4\,\ge\,m(\vp(x,t))\,\ge\,c_0\,\kappa(1-\kappa)\,>\,0 \quad\mbox{for every }\,
 (x,t)\in \overline{Q}.
 \end{equation}
This means that under the given hypotheses neither the possible degeneracy of $\,m\,$ nor the 
possible singularity of $\,f'\,$ 
can become active. Also, we may without loss of generality assume (by possibly choosing a larger
$K_1^*$) that
\begin{equation}
\label{mbounds}
\|m(\vp)\|_{C^0(\overline{Q})}\,+\,\|m'(\vp)\|_{C^0(\overline{Q})}\,
+\,\|m''(\vp)\|_{C^0(\overline{Q})}\,\le\,K^*_1\,.
\end{equation}    

\vspace{3mm}
{We will now show a global stability estimate. We have the following result.}

\vspace{3mm}
\noindent
{\sc \quad Proposition 2.2.}
\,\,{\em Let the hypotheses} {\bf (H1)}--{\bf (H8)} {\em be satisfied. Then there
exists a constant $\,K_2^*>0$, which only depends on the data of the 
state system and on $\,R$, such that it holds: whenever $\,\vb_1,\vb_2\in {\cal V}_R\,$ are given and} $\vp_1, \vp_2\in  C^1([0,T];L^2(\Omega))\cap H^1(0,T; H^1(\Omega))\cap 
L^\infty(0,T; H^2(\Omega))$ {\em  denote
 the associated solutions to the state system} (\ref{p1})--(\ref{p5}) {\em and $w_1\,,\,w_2$ the corresponding nonlocal operators according to} (\ref{p3}), {\em then we have for $\,\vp:=\vp_1-\vp_2$ and  $\vb:=\vb_1-\vb_2$, that 
for all $t\in [0,T]$ it holds
}
\begin{equation}
\label{stability}
\int_0^t\|\vp_t(s)\|_{L^2(\Omega)}^2\ds\,+\,\max_{0\le s\le t}\,
\|\vp(s)\|_{H^1(\Omega)}^2\,\leq \,K_2^*\,\int_0^t\|\vb(s)\|_{L^3(\Omega)^3}^2\ds \,.
\end{equation}

\proof In the following, the symbol $C$ will denote positive constants, which possibly differ from line to line or
even within lines. They may only depend on the problem data and $R$. 
To begin with, we put 
$$w:=w_1-w_2, \quad \mu_i:=f'(\varphi_i)+w_i, \,\,\,i=1,2, \quad
\mbox{and }\,\mu:=\mu_1-\mu_2,$$ 
and observe that $\,(\vp,w,\mu)\,$  
satisfies
\begin{align}
\label{j16}
&\vp_t\,-\,{\rm div}(m(\vp_1)\nabla\mu_1-m(\vp_2)\nabla\mu_2)\,=\,-\vb\cdot\nabla\vp_1\,-\,\vb_2\cdot\nabla\vp
\quad\,\mbox{a.e. in }\,Q\,,\\[2mm]
\label{j17}
&\mu=f'(\vp_1)-f'(\vp_2)+w\,\quad\mbox{a.e. in }\,Q\,,\\[2mm]
\label{j18}
&w(x,t)=-2\int_\Omega k(|x-y|)\vp(y,t)\dy \,\quad\mbox{a.e. in }\,Q\,,\\[2mm]
\label{j19}
&\vp(0)=w(0)=\mu(0)=0\,\quad\mbox{a.e. in }\,\Omega\,.
\end{align}
We also notice that, owing to (\ref{p2}), {\bf (H5)} and (\ref{flux}), we have 
\begin{eqnarray}
\label{j20}
&&\,\,\,(m(\vp_1)\nabla\mu_1-m(\vp_2)\nabla\mu_2)\, =\,c_0\,\nabla\vp\,+\,(m(\vp_1)\nabla w_1
-m(\vp_2)\nabla w_2)\\[2mm]
&&\,\,\,=c_0\,\nabla\vp \,+\,(m(\vp_1)-m(\vp_2))\nabla w_1 \,+\,m(\vp_2)\nabla w\,.\nonumber
\end{eqnarray} 
Hence, testing (\ref{j16}) by $\,\vp$, we have, for every $t>0$,
\begin{align}
\label{j21}
&\frac 12\,\|\vp(t)\|_{L^2(\Omega)}^2\,+\,c_0\int_0^t\!\!\!\int_\Omega|\nabla\vp|^2\dx\ds\,\le\,
\int_0^t\!\!\!\int_\Omega |\vb_2||\vp||\nabla\vp|\dx\ds\\[1mm]
&+\int_0^t\!\!\!\int_\Omega|m(\vp_1)-m(\vp_2)||\nabla w_1||\nabla\vp|\dx\ds
\,+\int_0^t\!\!\!\int_\Omega|m(\vp_2)||\nabla w||\nabla\vp|\dx\ds \nonumber\\[1mm]
&+\int_0^t\!\!\!\int_\Omega|\vb||\nabla\vp_1||\vp|\dx\ds\,.\nonumber
\end{align}

We denote the four integrals on the right-hand side by $I_j(t)$, $1\le j\le 4$, in that order, and estimate
them individually. At first, it follows from (\ref{defR}) and Young's inequality that
\begin{equation}
\label{j22}
I_1(t)\,\le\,\frac {c_0}8\int_0^t\!\!\!\int_\Omega|\nabla\vp|^2\dx\ds \,+\,
C\int_0^t\!\!\!\int_\Omega|\vp|^2\dx\ds\,.
\end{equation}
Next, from the mean value theorem, (\ref{mbounds}), (\ref{wbound}), and Young's inequality, we infer that
\begin{align}
\label{j23}
I_2(t)&\le \frac {c_0}8\int_0^t\!\!\!\int_\Omega|\nabla\vp|^2\dx\ds \,+\,C
\int_0^t\!\!\|\nabla w_1(s)\|^2_{L^\infty(Q)^3}\,\|\vp(s)\|^2_{L^2(\Omega)}\ds\\[1mm]
&\le \frac {c_0}8\int_0^t\!\!\!\int_\Omega|\nabla\vp|^2\dx\ds \,+\,
C\int_0^t\!\!\!\int_\Omega|\vp|^2\dx\ds\,. \nonumber
\end{align}
Moreover, (\ref{j1}) and Young's inequality imply that
\begin{align}
\label{j24}
I_3(t)&\le \frac {c_0}8\int_0^t\!\!\!\int_\Omega|\nabla\vp|^2\dx\ds\,+\,C\int_0^t\!\!\!\int_\Omega|\nabla w|^2\dx\ds
\\[1mm]
&\le \frac {c_0}8\int_0^t\!\!\!\int_\Omega|\nabla\vp|^2\dx\ds \,+\,C\int_0^t\!\!\!\int_\Omega|\vp|^2\dx\ds\,, \nonumber 
\end{align}
where the last inequality follows from {\bf (H7)}. Finally, we employ (\ref{phibound}), 
H\"older's and
Young's inequalities, as well as the continuity of the embedding $H^1(\Omega)\subset L^4(\Omega)$, to conclude that
\begin{align}
\label{j25}
I_4(t)&\le C\int_0^t \|\vb(s)\|_{L^2(\Omega)^3}\,\|\nabla\vp_1(s)\|_{L^4(\Omega)^3}\,\|\vp(s)\|_{L^4(\Omega)}\ds\\[1mm]
&\le \frac {c_0}8\int_0^t\|\vp(s)\|_{H^1(\Omega)}^2\ds\, +\, C \int_0^t \|\vp_1(s)\|_{H^2(\Omega)}^2
\,\|\vb(s)\|^2_{L^2(\Omega)^3}\ds\nonumber\\[1mm]
&\le \frac {c_0}8\int_0^t\!\!\!\int_\Omega|\nabla\vp|^2\dx\ds\, +\,C\int_0^t
\left(\|\vp(s)\|^2_{L^2(\Omega)}\,+\,\|\vb(s)\|^2_{L^2(\Omega)^3}\right)\ds\,.\nonumber
\end{align}
Combining the estimates (\ref{j21})--(\ref{j25}), and invoking Gronwall's lemma, we have thus shown that for any
$t\in [0,T]$ we have
\begin{equation}
\label{j26}
\max_{0\le s\le t}\,\|\vp(s)\|^2_{L^2(\Omega)}\,+\int_0^t\|\vp(s)\|^2_{H^1(\Omega)}\ds \,\le\,C
\int_0^t \|\vb(s)\|^2_{L^2(\Omega)^3}\ds\,,
\end{equation}
where the constant $C$ depends only on the data of the system and $R$.

Having established the stability estimate (\ref{j26}), we can now proceed to prove the stronger
estimate (\ref{stability}). To this end, we multiply (\ref{j16}) by $\vp_t$ and integrate over
$\Omega\times [0,t]$, where $t>0$. Integration by parts and (\ref{j20}) yield that  
\begin{align}\label{stability1}
&\int_0^t\!\!\!\int_\Omega|\vp_t|^2\dx\ds \,+\,\frac{c_0}2\,\|\nabla\vp(t)\|_{L^2(\Omega)^3}^2
\,\leq\,\sum_{i=1}^3 I_i(t)\,, 
\end{align} 
where
\begin{align}\no
&I_1(t)\,:=\,-\int_0^t\!\!\!\int_\Omega(m(\vp_1)\nabla w_1-m(\vp_2)\nabla w_2)\cdot \nabla\vp_t\dx\ds\,,\\[1mm] 
\no
&I_2(t)\,:=\,-\int_0^t\!\!\!\int_\Omega(\vb\cdot\nabla\vp_1)\,\vp_t\dx\ds\,,\\[1mm]
\no
&I_3(t)\,:=\,-\int_0^t\!\!\!\int_\Omega(\vb_2\cdot\nabla\vp)\,\vp_t\dx\ds\,.
\end{align}
We estimate these expressions individually. The last two terms are easily handled. Indeed,  
owing to (\ref{phibound}), H\"older's and Young's inequalities, and due to the continuity of
the embedding $H^1(\Omega)\subset L^6(\Omega)$, we have, for any $\gamma>0$ (to be specified later),
\begin{align}
\label{j27}
|I_2(t)|&\le \,C\int_0^t\|\vb(s)\|_{L^3(\Omega)^3}\,\|\nabla\vp_1(s)\|_{L^6(\Omega)^3}\,\|\vp_t(s)\|_{L^2(\Omega)}\ds\\[1mm]
&\le \,\gamma\int_0^t\!\!\!\int_\Omega|\vp_t|^2\dx\ds\,+\,\frac C\gamma \int_0^t
\|\vp_1(s)\|^2_{H^2(\Omega)}\,\|\vb(s)\|^2_{L^3(\Omega)^3}\ds \nonumber\\[1mm]
&\le \,\gamma\int_0^t\!\!\!\int_\Omega|\vp_t|^2\dx\ds\,+\,\frac C\gamma \int_0^t
\|\vb(s)\|^2_{L^3(\Omega)^3}\ds\,. \nonumber
\end{align}
Similarly, by also using {\eqref{defR} and} (\ref{j26}), we obtain that
\begin{align}
\label{j28}
|I_3(t)|&\le\int_0^t\|\vb_2(s)\|_{L^\infty(\Omega)^3}\,\|\nabla\vp(s)\|_{L^2(\Omega)^3}\,
\|\vp_t(s)\|_{L^2(\Omega)}\dx\ds \\[1mm]
&\le\, \,\gamma\int_0^t\!\!\!\int_\Omega|\vp_t|^2\dx\ds\,+\,\frac C\gamma {\int_0^t}
\|\vp{(s)}\|^2_{H^1(\Omega)}\ds \nonumber\\[1mm]
&\le \,\gamma\int_0^t\!\!\!\int_\Omega|\vp_t|^2\dx\ds\,+\,\frac C\gamma {\int_0^t}
\|\vb(s)\|^2_{L^2(\Omega)^3}\ds\,.\nonumber
\end{align}
It remains to estimate the first integral.
First notice that integration by parts with respect to time, together with \eqref{j19}, yields
\begin{align}
\no
I_1(t)&\,=\,-\int_\Omega \left((m(\vp_1)\nabla w_1-m(\vp_2)\nabla w_2)\cdot\nabla\vp\right)(t)\dx\\[1mm]
\no
&\,\,\,\quad+\int_0^t\!\!\!\int_\Omega(m(\vp_1)\nabla w_1-m(\vp_2)\nabla w_2)_t\cdot \nabla\vp\dx\ds
\,=:\,I_{11}(t)+I_{12}(t)\,.
\end{align}
Using the mean value theorem, (\ref{wbound}), (\ref{mbounds}), {\bf (H7)}, and Young's inequality, we obtain
\begin{align}
\label{j29}
|I_{11}(t)|&\le \int_\Omega |m(\vp_1(t))-m(\vp_2(t))||\nabla w_1(t)||\nabla\vp(t)|\dx\\[1mm]
&\quad + \int_\Omega |m(\vp_2(t))||\nabla w(t)||\nabla \vp(t)|\dx\nonumber\\[2mm]
&\le C\left(\|\varphi(t)\|_{L^2(\Omega)}\,\|\nabla w_1{(t)}\|_{L^\infty(\Omega)^3}
\,+\,\|\nabla w(t)\|_{L^2(\Omega)^3}\right)\|\nabla\vp(t)\|_{L^2(\Omega)^3}\nonumber\\[2mm]
&\le \,\gamma\,\|\nabla\vp(t)\|^2_{L^2(\Omega)^3}\,+\,\frac C \gamma \|\vp(t)\|^2_{L^2(\Omega)}
\nonumber\\[2mm]
&\le \,\gamma\,\|\nabla\vp(t)\|^2_{L^2(\Omega)^3}\,+\,\frac C \gamma \int_0^t\|\vb(s)\|^2_{L^2(\Omega)^3}\ds\,,\nonumber
\end{align}
where the last inequality follows from (\ref{j26}).

Finally, we estimate $I_{12}(t)$. We have
\begin{align}\label{j30}
|I_{12}(t)|&\,\leq \int_0^t\!\!\!\int_\Omega |m'(\vp_1)-m'(\vp_2)||\vp_{1,t}||\nabla w_1||\nabla\vp|\dx\ds\\[1mm]
&\quad +\int_0^t\!\!\!\int_\Omega |m'(\vp_2)||\vp_t||\nabla w_1||\nabla\vp|\dx\ds\nonumber\\[1mm]
&\quad +\int_0^t\!\!\!\int_\Omega |m'(\vp_2)||\vp_{2,t}||\nabla w||\nabla\vp|\dx\ds\nonumber\\[1mm]
&\quad +\int_0^t\!\!\!\int_\Omega |m(\vp_1)-m(\vp_2)||\nabla w_{1,t}||\nabla\vp|\dx\ds\nonumber\\[1mm]
&\quad +\int_0^t\!\!\!\int_\Omega |m(\vp_2)||\nabla w_t||\nabla\vp|\dx\ds\nonumber\\[1mm]
&\,=:\,J_1(t)+J_2(t)+J_3(t)+J_4(t)+J_5(t)\,.\nonumber
\end{align}
We estimate the terms on the right-hand side individually. First, invoking the mean value theorem,
H\"older's inequality, the continuity of the embedding $H^1(\Omega)\subset L^4(\Omega)$,
as well as (\ref{phibound}), (\ref{wbound}), (\ref{mbounds}) and (\ref{j26}), we find that 
for every $\gamma>0$ it holds
\begin{align}
\label{j31}
|J_1(t)|
&
\leq\,C\int_0^t \|\vp(s)\|_{L^4(\Omega)}\,\|\vp_{1,t}(s)\|_{L^4(\Omega)}\,\|\nabla w_1{(s)}\|_{L^\infty(\Omega)^3}\,\|\nabla\vp(s)\|_{L^2(\Omega)^3}\ds
\\[1mm]
&
\le\,C\,\max_{0\le s\le t}\,\|\vp(s)\|_{H^1(\Omega)}\Big(\int_0^t 
\|\nabla\vp(s)\|^2_{L^2(\Omega)^3}\ds\Big)^{1/2}\|\vp_{1,t}\|_{L^2(0,t;H^1(\Omega))}
\nonumber\\[1mm]
&
\le\,\gamma\,\max_{0\le s\le t}\,\|\vp(s)\|_{H^1(\Omega)}^2\,+\,\frac{C}{\gamma}
\int_0^t\|{\bf v}(s)\|^2_{L^2(\Omega)^3}\ds.
\nonumber
\end{align}
Moreover, using (\ref{wbound}), (\ref{mbounds}), (\ref{j26}) and Young's inequality, we have 
\begin{align}\label{j32}
|J_2(t)|&\le\,C\int_0^t \|m'(\vp_2(s))\|_{L^\infty(\Omega)}\,\|\vp_t{(s)}\|_{L^2(\Omega)}\,\|\nabla w_1(s)\|_{L^\infty(\Omega)^3}\,\|\nabla\vp(s)\|_{L^2(\Omega)^3}\ds\\[1mm]
&\leq \,\gamma \int_0^t\!\!\!\int_\Omega |\vp_t|^2\dx\ds\,+ \,\frac C 
\gamma\int_0^t \|\vb(s)\|^2_{L^2(\Omega)^3} \ds\,.\nonumber
\end{align}
Also, invoking (\ref{phibound}), (\ref{mbounds}), {\bf (H7)}, the continuity of the embedding $H^1(\Omega)\subset L^4(\Omega)$ and (\ref{j26}), we find the estimate
\begin{align}
\label{j33}
|J_3(t)|
&
\le\,\int_0^t\|m'(\vp_2(s)\|_{L^\infty(\Omega)}\,\|\vp_{2,t}(s)\|_{L^4(\Omega)}\, \|\nabla w(s)\|_{L^4(\Omega)^3}\,\|\nabla\vp(s)\|_{L^2(\Omega)^3}\ds\\[1mm]
&
\le\,C\int_0^t\|\vp_{2,t}(s)\|_{H^1(\Omega)}\,\|\vp(s)\|_{L^4(\Omega)}\,\|\nabla\vp(s)\|_{L^2(\Omega)^3}\ds 
\nonumber
\\
&
\le\,C\,\max_{0\le s\le t}\,\|\vp(s)\|_{H^1(\Omega)}\,\|\vp\|_{L^2(0,t;H^1(\Omega))}\,
\|\vp_{2,t}\|_{L^2(0,t;H^1(\Omega))}\nonumber
\\
&
\le\,\gamma \,\max_{0\le s\le t}\,\|\vp(s)\|_{H^1(\Omega)}^2\,+\,\frac C \gamma 
\int_0^t\|\vb(s)\|^2
_{L^2(\Omega)^3}\ds\,.\nonumber
\end{align}
Similarly, again using {\eqref{phibound}, \eqref{mbounds}, \eqref{j26}, together with} H\"older's inequality and {\bf (H7)}, we obtain that
\begin{align}
\label{j34}
|J_4(t)|&\le\,C\int_0^t\|\vp(s)\|_{L^4(\Omega)} \|\nabla w_{1,t}(s)\|_{L^4(\Omega)^3}\,\|\nabla \vp(s)\|_{L^2(\Omega)^3}\ds
\\[1mm]
&
\le\,C\int_0^t\|\vp(s)\|_{H^1(\Omega)}\,\|\vp_{1,t}(s)\|_{L^4(\Omega)}
\,\|\nabla\vp(s)\|_{L^2(\Omega)^3}\ds \nonumber
\\
&
\,\le\,\gamma \,\max_{0\le s\le t}\,\|\vp(s)\|_{H^1(\Omega)}^2\,+\,\frac C \gamma 
\int_0^t\|\vb(s)\|^2
_{L^2(\Omega)^3}\ds\,.\nonumber
\end{align}
Finally, we conclude from (\ref{mbounds}), {\bf (H7)}, {\eqref{j26},} as well as H\"older's and Young's inequalities, that
\begin{align}\label{j35}
|J_5(t)|&\le\,C\int_0^t\|\nabla w_t(s)\|_{L^2(\Omega)^3}\,\|\nabla\vp(s)\|_{L^2(\Omega)^3}\ds\\[1mm]
&\le\,C\int_0^t \|\vp_t(s)\|_{L^2(\Omega)}\,\|\vp(s)\|_{H^1(\Omega)}\ds\nonumber\\[1mm]
\nonumber
&\le\,\gamma\int_0^t\!\!\!\int_\Omega |\vp_t|^2\dx\ds\,+\,\frac C \gamma\int_0^t\|\vb(s)\|^2_{L^2(\Omega)^3}\ds\,.
\end{align}
Combining the estimates (\ref{stability1})--(\ref{j35}), 
and observing the continuity of the
embedding $L^3(\Omega)\subset L^2(\Omega)$, we have thus shown an estimate of the 
form
\begin{align}
&
(1-4\,\gamma)\int_0^t\!\!\!\int_\Omega |\vp_t|^2\dx\ds\,+\,\left(\frac{c_0}2\,-\gamma
\right)\|\nabla\vp(t)\|^2_{L^2(\Omega)^3}
\\
&
\quad\leq \,3\,\gamma\,\max_{0\le s\le t}
\,\|\vp(s)\|^2_{H^1(\Omega)}\,+\,\frac C \gamma \int_0^t\|\vb(s)\|^2_{L^3(\Omega)^3}\ds\,.\nonumber
\end{align}
From this, invoking (\ref{j26}), 
and adjusting $\gamma>0$ appropriately small,
it is easily seen that (\ref{stability}) is satisfied.
\qed

\vspace{5mm}
\noindent
{\bf Remark 3.} \,\,By virtue of {\bf (H7)} and {\bf (H8)}, the stability estimates (\ref{stability}) and
(\ref{j26}) entail corresponding estimates for $w$ and $\mu$. In particular, we may without loss of generality
assume (by possibly choosing an appropriately larger $K_2^*>0$) that for all $t\in [0,T]$ we have
\begin{align}
\label{stabw}
\int_0^t\|\nabla w_t(s)\|_{L^2(\Omega)^3}^2\ds\,+\,\|w\|_{L^\infty(0,T;H^2(\Omega))}^2&\,\le\,K_2^*\int_0^t\|\vb(s)\|
^2_{L^3(\Omega)^3}\ds\,,\\[1mm]  
\label{stabmu}
\int_0^t\|\nabla\mu_t(s)\|^2_{L^2(\Omega)^3}\ds\,+\,\|\mu\|^2_{L^\infty(0,T;H^2(\Omega))}&\,\le\,K_2^*\int_0^t\|\vb(s)\|
^2_{L^3(\Omega)^3}\ds\,.
\end{align}

\section{Optimal control}
\label{optimalcontrol}
\setcounter{equation}{0}

In this section, we study the optimal control problem {\bf (CP)} with ${\cal V}_{\rm ad}$ defined as in ${\bf (H1)}$, and
we assume that the general assumptions {\bf (H2)}--{\bf (H8)} are satisfied. Notice that, owing to Propositions 2.1 and 2.2, the {\em control-to-state operator} 
$${\cal S}: {\cal V}_R\to 
C^1([0,T];L^2(\Omega))\cap H^1(0,T;H^1(\Omega))\cap L^\infty(0,T;H^2(\Omega))
; \,\,\, \vb\mapsto\vp,
$$
is well defined and Lipschitz continuous as a mapping from $\mathcal{V}_R$ (viewed as a subset of $L^2(0,T; L^3(\Omega)^3)$) into $H^1(0,T; L^2(\Omega))\cap C^0([0,T]; H^1(\Omega))$. Moreover, all of the global bounds 
(\ref{phibound})--(\ref{mbounds}), as well as all of the stability estimates
(\ref{stability}), (\ref{j26}), (\ref{stabw}) and (\ref{stabmu}), are satisfied.

We are now ready to prove existence for the control problem {\bf (CP)}.
\bete\label{exicp}
Suppose that the hypotheses {\bf (H1)}--{\bf (H8)} are fulfilled. Then the problem {\bf (CP)} admits a solution $\bar \vb\in \mathcal{V}_{ad}$. 
\ente

\proof Let $\{\vb_n\}\subset\mathcal{V}_{ad}$ be a minimizing sequence for {\bf(CP)} and $\vp_n={\cal S}(\vb_n)$,
$n\in \NN$. Then it follows from {\bf (H1)} and (\ref{phibound}) that there exist $(\bar \vb, \bar\vp)$ such that, 
possibly for a subsequence which is again indexed by $n$, we have
\begin{eqnarray*}
\vb_n\to \bar \vb&&\quad\mbox{weakly in }\,L^2(0,T; H^1_{div}(\Omega)\cap H^1(0,T; L^3(\Omega)^3) \,,\\
&&\quad\mbox{and weakly-star in }\,L^\infty(Q)\\[1mm]
\vp_n\to\bar \vp &&\quad\mbox{weakly in }\, H^1(0,T; H^1(\Omega))\\
&&\quad\mbox{and weakly-star in }\, L^\infty(0,T; H^{2}(\Omega))\,,\\[1mm]
\partial_t\vp_n \to \partial_t\bar \vp&&\quad\mbox{weakly-star in }\,
L^\infty(0,T;L^2(\Omega)). 
\end{eqnarray*}

\noindent
Clearly, $\bar\vb\in {\cal V}_{\rm ad}$. In addition,
by virtue of standard compactness lemmas (cf. \cite[Thm. 5.1, p. 58]{lions} and \cite[Sec. 8, Cor. 4]{simon}), we have
the strong convergences 
\begin{align}\no
&\vb_n\to \bar \vb\quad\hbox{strongly in } C^0([0,T]; L^2(\Omega)^3)\,,\\\no
&\vp_n\to\bar\vp \quad\hbox{strongly in } C^0([0,T]; H^s(\Omega))\quad\forall s\in [0,2)\,, 
\end{align}
which implies, in particular, that
\begin{align}\no
&\vp_n\to\bar\vp \quad\hbox{strongly in } C^0(\overline{Q})\,,
\end{align}
as well as 
$$\vb_n\cdot\nabla\vp_n \to \bar\vb\cdot\nabla\bar\vp \quad\mbox{strongly in }\,L^1(Q).$$
Owing to the separation property (\ref{juerg2}) and the assumptions on $f$  and $m$, we also have
\begin{align}\no
&f'(\vp_n)\to f'(\bar\vp) \quad\mbox{and }\,\,\,m(\vp_n)\to m(\bar\vp), \quad\mbox{both strongly in }\,
C^0(\overline{Q})\,.
\end{align}
Finally, it is easily deduced from {\bf (H6)} that $\{w_n:={\cal K}(\vp_n)\}$ converges strongly in $C^0(\overline{Q})$ to $\bar w:={\cal K}(\bar\vp)$ (recall (\ref{p3}) and 
(\ref{intop})). 
In summary, we can pass to the limit as $n\to\infty$ in (\ref{p1})--(\ref{p5}), written for $(\vb_n,\vp_n)$, 
finding that $\bar\vp={\cal S}(\bar \vb)$; i.e., the pair $(\bar \vb, \bar \vp)$ is admissible for {\bf (CP)}. 
It then follows from the weak sequential lower semicontinuity properties of $J$ that 
$\bar\vb$, together with 
the associated state $\bar\vp={\cal S}(\bar \vb)$, is a solution to {\bf (CP)}. \qed 

\vspace{5mm}
We now turn our interest to the derivation of necessary first-order optimality conditions for problem {\bf (CP)}. 
Referring to \cite{tro} for a detailed discussion and description of the various techniques related to 
optimality conditions, we proceed as follows: we first prove a suitable differentiability property for the control-to-state 
operator ${\cal S}$, using the {\em linearized} system, and then we establish the necessary optimality
conditions in terms of a {\em variational inequality} and the associated {\em adjoint state equation}. In
the following, we will always (unless it is explicitly stated otherwise) assume that $\,\bar\vb\in {\cal V}_R\,$ is
fixed and that $\,(\bar\vp,\bar w,\bar \mu)\,$ is the associated triple solving the state system,
i.e.,
$$\bar\vp={\cal S}(\bar\vb), \quad\bar w={\cal K}(\bar\vp), \quad\bar\mu=f'(\bar\vp)+\bar w.$$ 

\paragraph{The linearized system} \,\,Suppose that an arbitrary $\,
{\bf h}\in {\cal V}\,$ is given. 
As a preparation for  the proof of  differentiability, we consider the following system, which 
is obtained by linearizing the state system (\ref{p1})--(\ref{p5}) at $\bar\vp=S(\bar \vb)$:
\begin{eqnarray}\label{l1}
&&\xi_t\,-\,c_0\,\Delta\xi-\dive\left(m'(\bar\vp)\,\xi\,\nabla\bar w
\,-2\,m(\bar\vp)\,\nabla\left(\int_\Omega k(|x-y|)\,\xi(y,\,\cdot\,)\dy\right)\right)\\[2mm]
&&=-\,{\bf h}\cdot \nabla \bar\vp\,-\,\bar \vb\cdot\nabla\xi\,\quad\hbox{a.e. in }Q\,,\nonumber\\[2mm]
\label{l2}
&&\bar w(x,t)=\int_\Omega k(|x-y|)(1-2\bar\vp(y,t))\dy\quad\hbox{a.e. in }Q\,,\\[2mm]
\label{l3}
&&\Big(c_0\,\nabla\xi+m'(\bar\vp)\,\xi\,\nabla\bar w \,-2\,m(\bar\vp)\,\nabla\Big(\int_\Omega k(|x-y|)\,\xi(y,\,\cdot\,)\,dy\Big)\Big)\cdot{\bf n}=0\\[1mm]
&&\quad\,\hbox{ a.e. on }\,\Sigma\,,\nonumber\\[2mm]
\label{l4}
&&\xi(0)=0\quad\hbox{a.e. in }\Omega\,.
\end{eqnarray}

\noindent
After proving that (\ref{l1})--(\ref{l4}) has a unique solution $\xi$, we expect that 
$\,\xi=D{\cal  S}(\bar \vb){\bf h}$, where $\,D{\cal S}(\bar\vb)\,$ denotes the Fr\'echet
derivative of $\,{\cal S}\,$ at $\,\bar\vb$. Recalling the global bounds (\ref{phibound})--(\ref{mbounds}), we 
can expect the regularity 
\begin{equation}\label{regoxi}
\xi\in H^1(0,T; L^2(\Omega))\cap C^0([0,T]; H^1(\Omega))\cap L^2(0,T; H^2(\Omega))\,.
\end{equation}
We have the following result.

\vspace{5mm}
\bepr\label{prop-exilin}
Suppose that the hypotheses {\bf (H1)}--{\bf (H8)} are fulfilled. Then the problem 
{\rm (\ref{l1})}--{\rm (\ref{l4})} has a unique solution satisfying {\rm (\ref{regoxi})}. 
\empr
\proof The proof is performed via a Faedo--Galerkin scheme inspired by \cite[Sec. 4]{cfg}. 
To  this end, we choose $\,\{\psi_j\}_{j\in\NN}\,$ to be the 
family of (appropriately orthonormalized and ordered) eigenfunctions to the eigenvalue problem
$$-\Delta\psi+\psi=\lambda\,\psi\,\quad\mbox{in }\,\Omega\,,\qquad \frac{\partial\psi}
{\partial {\bf n}}=0\,\quad\mbox{on }\,\partial\Omega\,,$$ 
as a Galerkin basis in $H^1(\Omega)$. 
Putting $\,\xi_n(x,t):=\sum_{k=1}^{n} a_k (t)\psi_k(x)$, we then
look for a solution to the approximating problem 
\begin{align}\label{p1n}
&\int_\Omega \xi_n'(t)\,\psi\dx+\int_\Omega c_0\nabla\xi_n(t)\cdot\nabla\psi \dx+\int_\Omega\Big(m'(\bar\vp(t))\,\xi_n(t)\,\nabla\bar w(t)\\[1mm]
\no
&-\,2\,m(\bar\vp(t))\,\nabla
\Bigl(\int_\Omega k(|x-y|)\xi_n(y,t)\dy\Bigr)\Bigr)\cdot\nabla\psi\dx\quad\\[1mm]
\no
&\, =-\int_\Omega ({\bf h}(t)\cdot \nabla \bar\vp(t))\,\psi\dx-\int_\Omega(
\bar \vb(t)\cdot\nabla\xi_n(t))\,\psi\dx\,,\quad\,\mbox{for }\,t\in (0,T], \\[2mm]
\label{p2n}
&\qquad \qquad \xi_n(0)=0\quad \hbox{a.e. in }\Omega\,,
\end{align}
for every $\psi\in \Psi_n:=\mbox{span}\,\{\psi_1, \dots, \psi_n\}$. Apparently, this
is nothing but an initial value problem for a system of linear ordinary differential equations for the unknown functions
$\,a_1,...,a_n$, where, owing to the global bounds (\ref{phibound})--(\ref{mbounds}), all occurring coefficient
functions are continuous on $[0,T]$. It is therefore a standard matter to show   
that there exists some $\,T_n\in (0,T]\,$ such that the ODE system has a maximal solution ${\bf a}:=(a_1, \dots, a_n)\in C^1([0,T_n); \RR^n)$ that specifies a solution $\,\xi_n\in C^1([0,T_n);H^3(\Omega))$. Observe that
\begin{equation}
\label{bcxi}
\frac{\partial \xi_n}{\partial {\bf n}} =\frac{\partial \Delta\xi_n}{\partial {\bf n}}=0\,\quad
\mbox{on }\,\Sigma \quad\forall\,n\in\NN\,.
\end{equation} 

We now aim to prove a (uniform in $n\in\NN$) estimate for $\xi_n$ in  $C^0([0,T]; H^1(\Omega))\cap L^2(0,T; H^2(\Omega))$. 
Once this is shown, it is a standard matter to show that $T_n=T$ and to pass to the limit as $n\to\infty$ 
to recover a solution with the
asserted regularity to the linearized problem (\ref{l1})--(\ref{l4}). Due to regularity of the coefficients and of the known functions in the system, we can also prove that the solution is unique, simply by testing the difference between two equations (\ref{l1}), written for two possible different solutions $\xi_1$ and $\xi_2$, by $\xi_1-\xi_2$ and then 
exploiting the linearity of the problem. Since
these arguments are straightforward, we can 
allow ourselves to be brief here and to restrict ourselves 
to the derivation of the asserted global bounds. 

To this end, let $t\in (0,T_n)$ be arbitrary. In what follows, $C_i$, $i\in\NN$, will denote positive constants that
may depend on the data of the system but not on $n\in\NN$, and we will make repeated use of
the global bounds (\ref{phibound})--(\ref{mbounds}) {and of \eqref{defR}}. First observe that the global  bounds and hypothesis 
{\bf(H7)} imply that, for any $s\in [0,t]$,
\begin{eqnarray}
\label{jneu}
&&\Big \|m'(\bar\vp(s))\,\xi_n(s)\,\nabla\bar w(s) \,-\,2\,m(\bar\vp(s))\,\nabla
\Bigl(\int_\Omega k(|x-y|)\xi_n(y,s)\dy\Bigr)\Big\|_{L^2(\Omega)}\\[2mm]
&&\le\,C_1\,\|\xi_n(s)\|_{L^2(\Omega)}\,.\nonumber
\end{eqnarray}  

Now we insert $\,\psi=\xi_n(t)\,$ in (\ref{p1n}) and integrate over $[0,t]$ to find that 
\begin{equation}
\label{j40}
\frac 12\,\|\xi_n(t)\|^2_{L^2(\Omega)}\,+\,c_0\int_0^t\!\!\int_\Omega|\nabla\xi_n|^2\dx\ds
\,\le\,\sum_{j=1}^3 I_j(t)\,,
\end{equation}
with expressions $I_j(t)$, $1\le j\le 3$, that will be specified and estimated below. 

Let $\gamma>0$ be arbitrary (to be specified later). We have, using (\ref{jneu}) and
Young's inequality,
\begin{eqnarray}
\label{bI1}
&&\\[-1mm]
&&|I_1(t)|=\Big |\!\int_0^t\!\!{\int_\Omega\!\!} \left[m'(\bar\vp)\,\xi_n\,
\nabla\bar w -2\,m(\bar\vp)\,\nabla
\Bigl(\int_\Omega\!\! k(|x-y|)\xi_n(y,s)\dy\Bigr)\right]\cdot\nabla\xi_n\dx\ds
\Big |\nonumber\\[1mm]
&&\hspace*{12.5mm}\le\,\gamma \int_0^t\!\!\int_\Omega|\nabla \xi_n|^2\dx\ds\,+\,\frac{C_2}\gamma
\int_0^t\!\!\int_\Omega|\xi_n|^2\dx\ds\,.\nonumber
\end{eqnarray}

\noindent Moreover, 
\begin{eqnarray}
\label{bI3}
&&\\[-1mm]
&&|I_2(t)|\,=\,\Big|\int_0^t\!\!\int_\Omega \!\!({\bf h}\cdot\nabla\bar\vp)\,\xi_n\dx\ds\Big|
\,\le\,\int_0^t\!\!\|{\bf h}(s)\|_{L^3(\Omega)^3}\,\|\nabla\bar\vp(s)\|_{L^6(\Omega)^3}
\,\|\xi_n(s)\|_{L^2(\Omega)}\ds\nonumber\\[1mm]
&&\hspace*{12.5mm} \le\,C_3\,\Big(\int_0^t\|{\bf h}(s)\|_{L^3(\Omega)^3}^2\ds
\,+\,\int_0^t\!\!\int_\Omega|\xi_n|^2\dx\ds\Big)\,,\nonumber
\end{eqnarray}
as well as
\begin{eqnarray}
\label{bI4}
&&|I_3(t)|\,=\,\Big|\int_0^t\!\!\int_\Omega (\bar\vb\cdot\nabla\xi_n)\,\xi_n\dx\ds\Big|
\\[1mm]
&&\hspace*{12.5mm}\le \,C_4\int_0^t \|\bar\vb(s)\|_{L^\infty(\Omega)^3}\,\|\nabla\xi_n(s)\|_{L^2
(\Omega)^3}\,\|\xi_n(s)\|_{L^2(\Omega)}\ds\nonumber\\[1mm]
&&\hspace*{12.5mm} \le\,\gamma\int_0^t\!\!\int_\Omega|\nabla\xi_n|^2\dx\ds
\,+\,\frac{C_5}\gamma\int_0^t\!\!\int_\Omega|\xi_n|^2\dx\ds\,.\nonumber
\end{eqnarray}
Combining the estimates (\ref{j40})--(\ref{bI4}), choosing $\gamma>0$ small
enough, and applying Gronwall's lemma, we have thus shown the estimate
\begin{equation}
\label{j41}
\max_{0\le s\le t}\,\|\xi_n(s)\|^2_{L^2(\Omega)} \,+\int_0^t 
\|\xi_n(s)\|^2_{H^1(\Omega)}\ds\,\le\,C_{6}\int_0^t\|{\bf h}(s)\|^2_{L^3(\Omega)^3}\ds
\,\le\,C_{7}\,.
\end{equation}

Next, we insert $\psi=-\Delta\xi_n(t)$ in (\ref{p1n}), integrate by parts using
the boundary condition (\ref{bcxi}), and then integrate over $[0,t]$. We then obtain
\begin{eqnarray}
\label{j42}
\qquad\qquad\frac 12 \|\nabla\xi_n(t)\|^2_{L^2(\Omega)^3}\,+\,c_0\int_0^t\!\!\int_\Omega
|\Delta\xi_n|^2\dx\ds\,\le\,\int_0^t\!\!\int_\Omega(|g_1|\,+\,|g_2|)\,|\Delta\xi_n|\dx\ds\\[1mm]
+\int_0^t\!\!\int_\Omega |{\bf h}||\nabla\bar\vp||\Delta\xi_n|\dx\ds \,+\int_0^t\!\!\int_\Omega |\bar\vb||\nabla\xi_n||\Delta\xi_n|\dx\ds \,,\nonumber\hspace*{20mm}
\end{eqnarray}
where the functions $g_1$, $g_2$ will be specified and estimated below. Now let $\gamma>0$ be arbitrary 
(to be specified later). The last two {integrals} on the right-hand side of (\ref{j42}) are easily 
estimated. In fact, using the general bounds (\ref{phibound}), as well as H\"older's and Young's
inequalities {and \eqref{defR}}, we have
\begin{eqnarray}
\label{j43}
\\[-1mm]
\int_0^t\!\!\int_\Omega |{\bf h}||\nabla\bar\vp||\Delta\xi_n|\dx\ds \,\le\,
C_{8}\int_0^t \|{\bf h}(s)\|_{L^3(\Omega)^3}\,\|\nabla\bar\vp(s)\|_{L^6(\Omega)^3}\,
\|\Delta\xi_n(s)\|_{L^2(\Omega)}\ds\nonumber\\[1mm]
\le\,\gamma\int_0^t\!\!\int_\Omega |\Delta\xi_n|^2\dx\ds\,+\,\frac{C_{9}}\gamma\int_0^t\|{\bf h}(s)\|^2_{L^3(\Omega)^3}
{\ds} \,,\nonumber\hspace*{2mm} 
\end{eqnarray} 
as well as
\begin{eqnarray}
\label{j44}\\[-1mm]
\int_0^t\!\!\int_\Omega |{\bar\vb}||\nabla\xi_n||\Delta\xi_n|\dx\ds \,\le\,\int_0^t\|\bar\vb(s)\|_{{L^\infty(\Omega)^3}}
\,\|\nabla\xi_n(s)\|_{{L^2(\Omega)^3}}\,\|\Delta\xi_n(s)\|_{L^2(\Omega)}\ds\nonumber\\[1mm]
\le\,\gamma\int_0^t\!\!\int_\Omega |\Delta\xi_n|^2\dx\ds\,+\,\frac{C_{10}}\gamma\int_0^t\|\xi_n(s)\|^2_{H^1(\Omega)}
\ds\,.\nonumber
\end{eqnarray}

It remains to estimate the first integral on the right-hand side of (\ref{j42}). To this end, we first infer
from the global bounds (\ref{phibound})--(\ref{mbounds}) that a.e. on $Q$ it holds
\begin{eqnarray}
\label {j45}
|g_1|&\!\!:=\!\!&\left|\,{\rm div}\,\left[m'(\bar\vp)\,\xi_n\,\nabla\bar w
\right]\,\right|\,\le\,C_{11}\Bigl(|\xi_n|\,(|\nabla\bar\vp|\,+\,
|\Delta\bar w|)\,+\,|\nabla\xi_n||\nabla\bar\vp|\Bigr)\,,  \nonumber
\end{eqnarray} 
where it is easily verified that the expression in the inner bracket,
which we denote by $z$, is bounded in
$C^0([0,T];L^2(\Omega))$. We thus have, invoking (\ref{young}) and (\ref{compactness2}),
\begin{eqnarray}
\label{j46}
&&\int_0^t\!\!\int_\Omega |\xi_n||z||\Delta\xi_n|\dx\ds\,\le\,
\int_0^t\|\xi_n(s)\|_{L^\infty(\Omega)}\,\|z(s)\|_{L^2(\Omega)}\,\|\Delta\xi_n(s)\|_{L^2(\Omega)}\ds
\\[1mm]
&&\quad\le \,\gamma\int_0^t\!\!\int_\Omega |\Delta\xi_n|^2\dx\ds\,+\,\frac{C_{12}}
\gamma\int_0^t\|\xi_n(s)\|^2_{L^\infty(\Omega)}\ds\nonumber\\[1mm]
&&\quad\le \,2\,\gamma\int_0^t\|\xi_n(s)\|^2_{H^2(\Omega)}\ds\,+\,C_{13}\,(\gamma^{-1}+\gamma^{-3})
\,\int_0^t\|\xi_n(s)\|^2_{H^1(\Omega)}\ds\,.
\nonumber 
\end{eqnarray}  
Moreover, by (\ref{phibound}), (\ref{compactness}), and H\"older's and Young's inequalities,
it holds
\begin{eqnarray}
\label{j47}
\\[-1mm]
\int_0^t\!\!\int_\Omega \!\!|\nabla\bar\vp||\nabla\xi_n||\Delta\xi_n|\dx\ds\,\le\,C_{14}\!\!
\int_0^t\!\!\|\nabla\bar\vp(s)\|_{{L^6(\Omega)^3}}\,\|\nabla\xi_n(s)\|_{L^3(\Omega)^3}\,\|\Delta\xi_n(s)\|
_{L^2(\Omega)}\ds\nonumber\\[1mm]
\le\, \gamma\int_0^t\!\!\int_\Omega |\Delta\xi_n|^2\dx\ds\,+\,\frac{C_{15}}
\gamma\int_0^t\|\nabla\xi_n(s)\|^2_{L^3(\Omega)^3}\ds \hspace*{40mm}\nonumber\\[1mm]
\le\,2\,\gamma\int_0^t\|\xi_n(s)\|^2_{H^2(\Omega)}\ds\,+\,C_{16}\,(\gamma^{-1}+\gamma^{-3})
\,\int_0^t\|\xi_n(s)\|^2_{H^1(\Omega)}\ds\,.\hspace*{16.5mm}\nonumber
\end{eqnarray}
Finally, notice that {for a.e. $(x,t)\in  Q$} we have
\begin{eqnarray*}
|g_2(x,t)|\,:=\,\left|\,{\rm div}\,\left[2\,m(\bar\vp(x,t))\,\nabla{\left(\int_\Omega
k(|x-y|)\,\xi_n(y,t)\dy\right)}\right]\,\right|\,,
\end{eqnarray*}
and it easily follows from (\ref{phibound}), (\ref{mbounds}), and the hypotheses
{\bf (H7)} and {\bf (H8)}, that 
$$\int_0^t\|g_2(s)\|^2_{L^2(\Omega)}\ds\,\le\,C_{17}\int_0^t\|\xi_n(s)\|^2_{H^1(\Omega)}\ds\,,$$
whence we obtain that
\begin{equation}
\label{j48}
\int_0^t\!\!\int_\Omega|g_2||\Delta\xi_n|\dx\ds\,\le\,\
\gamma\int_0^t\!\!{\int_\Omega}|\Delta\xi_n|^2\dx\ds\,+\,\frac{C_{18}}\gamma\int_0^t\|\xi_n(s)\|^2_{H^1(\Omega)}\ds\,.
\end{equation} 
Now observe that $\partial\xi_n/\partial{\bf n}=0$, so that standard elliptic estimates imply
that
$$\|\xi_n(s)\|_{H^2(\Omega)}\,\le\,C_{19}\left(\|\Delta\xi_n(s)\|_{L^2(\Omega)}\,+\,
\|\xi_n(s)\|_{H^1(\Omega)}\right)
\,,$$
where $C_{19}>0$ depends only on $\Omega$. Therefore, choosing $\gamma>0$ appropriately small,
and invoking (\ref{j41}), we can infer from the estimates (\ref{j42})--(\ref{j48}) that
\begin{equation}
\label{j49}
\max_{0\le s\le t}\|\xi_n(s)\|_{H^1(\Omega)}^2\,+\int_0^t\|\xi_n(s)\|^2_{H^2(\Omega)}\ds
\,\le\,C_{20}\int_0^t\|{\bf h}(s)\|^2_{L^3(\Omega)^3}\ds\,\le\,C_{21}\,.
\end{equation}
This concludes the proof of the assertion.\qed

\vspace{7mm}
\noindent
{\bf Remark 4.} \,\,From (\ref{j49}) it follows, in particular, that the linear mapping 
$\,{\bf h}\mapsto \xi=:\xi^{\bf h}$\, is continuous as a mapping from ${\cal V}$ into
the space $\,C^0([0,T];H^1(\Omega))\cap L^2(0,T;H^2(\Omega))$. 

\vspace{5mm}
\paragraph{Differentiability of the control-to-state mapping} In this section we are going
to prove the following result:

\vspace{5mm}
\bepr\label{prop-differ}
Let the hypotheses {\bf (H1)}--{\bf (H8)} be satisfied. Then the control-to-state operator 
$${\cal S}\,:\, {\cal V}_R\to 
C^1([0,T];L^2(\Omega))\cap H^1(0,T;H^1(\Omega))\cap L^\infty(0,T;H^2(\Omega)),
\quad {\bf v}\mapsto \vp\,,$$ 
is Fr\'echet differentiable in ${\cal V}_R$ as a mapping from ${\cal V}$ into ${\cal Y}:= 
C^0([0,T];L^2(\Omega))\cap L^2(0,T; H^1(\Omega))$, and for every $\bar\vb\in {\cal V}_R$ the Fr\'echet derivative $D{\cal S}(\bar\vb)\in {\cal L}({\cal V},
{\cal Y})$ is defined as follows: for every $\,{\bf h}\in {\cal V}\,$ we have
\begin{equation}
\label{Frechet}
D{\cal S}(\bar\vb){\bf h}=\xi^{\bf h}\,,
\end{equation}
where $\,\xi^{\bf h}\,$ is the unique solution to the linearized system 
{\rm (\ref{l1})--(\ref{l4})} with
$\bar\vp={\cal S}(\bar\vb)$. 
\empr
\proof Let $\,\bar\vb\in {\cal V}_R\,$ be fixed, and let $\bar\vp={\cal S}(\bar\vb)$. Since ${\cal V}_R$ is open, there is some $\Lambda>0$
such that $\bar\vb+{\bf h}\in {\cal V}_R$ whenever $\,\|{\bf h}\|_{{\cal V}}\le\Lambda$. In the following,
we only consider such perturbations ${\bf h}$ and set 
\[
\vb^{\bf h} =\bar \vb+{\bf h}, \quad \vp^{\bf h}={\cal S}(\vb^{\bf h}), \quad y^{\bf h}=\vp^{\bf h}-\bar\vp-\xi^{\bf h}\,. 
\]
Since the linear mapping $\,{\bf h}\mapsto \xi^{\bf h}$\, is by Remark 4 continuous as a mapping from
${\cal V}$ into ${\cal Y}$, it suffices to show that 
there exists an increasing mapping  $Z\,:\, (0,+\infty)\to (0,+\infty)$ such that $\,\lim_{\lambda\searrow0}Z(\lambda)/\lambda^2=0$ and 
\begin{equation}\label{diffe}
\|\yl\|^2_{ C^0([0,T]; L^2(\Omega))\cap L^2(0,T; H^1(\Omega))}\leq Z(\|{\bf h}\|_{\mathcal{V}})\,. 
\end{equation}

In the following, we will denote by $C_i$, $i\in\NN$, positive constants that may
depend on the data of the system and on $R$, but not on the special choice of
${\bf h}\in {\cal V}$ with $\,\|{\bf h}\|_{\cal V}\le\Lambda$. For a shorter exposition, we also
often omit the arguments of the involved functions if no confusion may arise. Notice that the global 
bounds (\ref{phibound}), (\ref{f'bounds}) and
(\ref{mbounds}) are satisfied by $\,\vp^{\bf h}\,$ for any perturbation ${\bf h}$ with $\|{\bf h}\|_{\cal V}
\le\Lambda$, and, owing to the weak sequential lower semicontinuity of norms, it follows from (\ref{j49})
that for all such perturbations we have 
\begin{equation}
\label{j50}
\|\xi^{\bf h}\|_{C^0([0,T];H^1(\Omega))\cap L^2(0,T;H^2(\Omega))}\,\le\,C_1\,.
\end{equation}   
 
First of all, it is easily verified that $\yl$ is a 
strong solution to the following system:
\begin{align}\label{y1}
&\yl_t-c_0\,\Delta\yl -\dive\left[m(\vpl)\nabla\Big(\int_\Omega k(|x-y|)(1-2\vpl)\dy\Big)\right.\\
\no
&\hspace*{35mm} -m(\bar \vp)\nabla\Big(\int_\Omega k(|x-y|)(1-2\bar\vp)\dy\Big)\\
\no
&\hspace*{35mm} -m'(\bar\vp)\,\xi^{\bf h}\,\nabla\Big(\int_\Omega k(|x-y|)(1-2\bar\vp)\dy\Big)\\
\no
&\hspace*{35mm}\left.
+ \,2\,m(\bar\vp)\nabla\Big(\int_\Omega k(|x-y|)\,\xi^{\bf h}(y)\dy\Big)\right]\\[2mm]
\no
&+\,\vb\cdot\nabla\yl\,+\,{\bf h}\cdot (\nabla \vpl-\nabla\bar\vp)=0\quad\hbox{a.e. in }Q\,,
\end{align}
\begin{align}
\label{y2}
&\left[c_0\,\nabla\yl \,+\,m(\vpl)\nabla\Big(\int_\Omega k(|x-y|)(1-2\vpl)\dy\Big)\right.\\
\no
&\quad -m(\bar \vp)\nabla\Big(\int_\Omega k(|x-y|)(1-2\bar\vp)\dy\Big)\\
\no
&\quad -m'(\bar\vp)\,\xi^{\bf h}\,\nabla\Big(\int_\Omega k(|x-y|)(1-2\bar\vp)\dy\Big)\\
\no
&\quad\left.
+2\,m(\bar\vp)\nabla\Big(\int_\Omega k(|x-y|)\xi^{\bf h}(y)\dy\Big)
\right]\cdot{\bf n}=0\quad\hbox{a.e. on }\Sigma\,,\\[2mm]
\label{y3}
&\yl(0)=0\quad\hbox{a.e. in }\Omega\,.  
\end{align}
We now test \eqref{y1} by $\yl$, integrate over $(0,t)$ where $t\in (0,T]$, and use \eqref{y2} and \eqref{y3} to get 
\begin{eqnarray}
\label{j51}
&&\frac12\|\yl(t)\|_{L^2(\Omega)}^2\,+\int_0^t\!\!\int_\Omega |\nabla\yl|^2\dx\ds\\
&&+\int_0^t\!\!\int_\Omega \nabla \yl\cdot \left\{m(\vpl)\nabla\Big(\int_\Omega k(|x-y|)(1-2\vpl)\dy\Big)\right.
\nonumber\\
&&\hspace*{27mm} -m(\bar \vp)\nabla\Big(\int_\Omega k(|x-y|)(1-2\bar\vp)\dy\Big)\nonumber\\
&&\hspace*{27mm} -m'(\bar\vp)\,\xi^{\bf h}\,\nabla\Big(\int_\Omega k(|x-y|)(1-2\bar\vp)\dy\Big)\nonumber\\
&&\hspace*{27mm}\left.
+\,2\,m(\bar\vp)\,\nabla\Big(\int_\Omega k(|x-y|)\xi^{\bf h}(y)\dy\Big)
\right\}\dx\ds\nonumber\\
&&+\int_0^t\!\! \int_\Omega\yl(\bar \vb\cdot\nabla \yl)\dx\ds\,+
\int_0^t\!\!\int_\Omega {\bf h}\cdot(\nabla \vpl-\nabla\bar\vp)\, \yl\dx\ds\,=\,0\,.\nonumber
\end{eqnarray} 
We have $\,\int_0^t\!\! \int_\Omega\yl\,(\bar \vb\cdot\nabla \yl)\dx\ds=0$, since
$\bar\vb$ vanishes on $\partial\Omega$ and is divergence free. Moreover, using H\"older's 
and Young's inequalities, as well as
the stability estimate (\ref{stability}) and the continuity of the embedding $H^1(\Omega)\subset
L^6(\Omega)$, we obtain that 
\begin{align}
\label{i1}
&\left|\int_0^t\!\!\int_\Omega {\bf h}\cdot(\nabla \vpl-\nabla\bar\vp) \yl\dx\ds\right|\,\leq
 \int_0^t\!\!\|\nabla (\vpl-\bar\vp)(s)\|_{{L^2(\Omega)^3}}\|\yl(s)\|_{L^6(\Omega)}\|{\bf h}(s)\|_{{L^3(\Omega)^3}}\ds\\
\no
&\qquad\leq \gamma \int_0^t\|\yl(s)\|_{H^1(\Omega)}^2\ds
\,+\,\frac{C_2}\gamma\int_0^t\|{\bf h}(s)\|_{L^3(\Omega)}^2\,\|\nabla(\vpl-\bar\vp)(s)\|_{L^2(\Omega)}^2\ds\\
\no
&\qquad\leq \gamma \int_0^t\|\yl(s)\|_{H^1(\Omega)}^2\ds\,+\,
\frac{C_3}\gamma\Big(\int_0^t\|{{\bf h}}(s)\|_{L^3(\Omega)}^2\ds \Big)^2\,,
\end{align}
for every positive $\gamma$ (to be chosen later). 

It remains to estimate the third summand in (\ref{j51}). To this end, we observe that 
the expression in the curly bracket in (\ref{j51}) equals the sum of
the following 
three expressions:
\begin{eqnarray*}
&&A_1(x,s):=\left(m(\vp^{\bf h})-m(\bar\vp)-m'(\bar\vp)\,\xi^{\bf h}\right)\!(x,s)\,\nabla\!
\int_\Omega\! k(|x-y|)\,(1-2\,\bar\vp(y,s))\dy, \\[1mm]
&&A_2(x,s):=\,-2\,(m(\vp^{\bf h})-m(\bar\vp))(x,s)\,\nabla\!\int_\Omega\!k(|x-y|)(\vp^{h}(y,s)-\bar\vp(y,s))\dy,\\[1mm]
&&A_3(x,s):=\,-2\,m(\bar\vp(x,s))\,\nabla\!\int_\Omega \!k(|x-y|)\,y^{\bf h}(y,s)\dy\,.
\end{eqnarray*}
Moreover, Taylor's theorem, using also the separation property (\ref{juerg2}) and the global bounds (\ref{mbounds}), yields
that almost everywhere in  $Q$ it holds
\begin{align}\label{taylor}
|m(\vp^{\bf h})-m(\bar\vp)-m'(\bar\vp)\xi^{\bf h}|&\le\,C_4\,|y^{\bf h}|\,+\,\frac 12\,\max_{\kappa\le\sigma\le 1-\kappa}
\,|m''(\sigma)||\vp^{\bf h}-\bar\vp|^2\\[1mm]
&\le\,C_4\,|y^{\bf h}|\,+\,C_5\,|\vp^{\bf h}-\bar\vp|^2\,.\nonumber
\end{align}

\noindent
Now, by virtue of hypothesis {\bf (H7)} and (\ref{wbound}), and by invoking H\"older's and Young's inequalities, we have
\begin{align}\label{i21}
&\int_0^t\!\!\int_\Omega |A_1||\nabla y^{\bf h} |\dx\ds\,\le\,C_6\int_0^t{\int_\Omega}\left(|y^{\bf h}|\,+\,|\vp^{\bf h}-\bar\vp|^2\right)|\nabla \yl|\dx\ds\\[2mm]
\no
&\,\leq \int_0^t\!\!\int_\Omega\!\!
\Big(\gamma|\nabla\yl|^2\,+\,\frac{C_7}\gamma\,|\yl|^2\Big)\dx\ds \,+\int_0^t\Big(\|\vpl-\bar\vp\|_{L^4(\Omega)}^2\|\nabla\yl\|_{{L^2(\Omega)^3}}\Big)(s)\ds\\[2mm]
\no
&\,\leq 2\,\gamma\int_0^t\!\!\int_\Omega\!|\nabla\yl|^2\dx\ds \,+\,\frac{C_8}\gamma
\Big(\int_0^t\!\!\int_\Omega\!|\yl|^2\dx\ds\,+\Big(\int_0^t\|{\bf h}(s)\|_{L^3(\Omega)^3}^2\ds\Big)^2\Big)\,,
\end{align}
where again (\ref{stability}) was employed. 
Similarly,
\begin{align}\label{i22}
&\int_0^t\!\!\int_\Omega |A_2||\nabla y^{\bf h} |\dx\ds\,\le\, C_{9}
\int_0^t\|\nabla y^{\bf h}(s)\|_{L^2(\Omega)^3}\,\|\vp^{\bf h}(s)-\bar\vp(s)\|^2_{L^4(\Omega)}\ds \\[2mm]
\no
&\le\,\gamma\int_0^t\!\!\int_\Omega |\nabla y^{\bf h}|^2\dx\ds\,+\,\frac{C_{10}}\gamma
\int_0^t\|\vp^{\bf h}(s)-\bar\vp(s)\|^4_{L^4(\Omega)}\ds\,,
\end{align}
as well as, using hypothesis {\bf (H7)} once more, 
\begin{align}\label{j52}
&\int_0^t\!\!\int_\Omega |A_3||\nabla y^{\bf h} |\dx\ds\,\le\,C_{11}\int_0^t
\|\nabla \yl(s)\|_{L^2(\Omega)^3}\,\Big\|\nabla\int_\Omega k(|x-\eta|)\yl(\eta,s)\,{\rm d}\eta\Big\|_{L^2(\Omega)^3}\ds\\[2mm]
\no
&\le\,\gamma\int_0^t\!\!\int_\Omega |\nabla y^{\bf h} |^2\dx\ds \,+\,  \frac {C_{12}}\gamma
\int_0^t\!\!\int_\Omega |y^{\bf h} |^2\dx\ds\,.
\end{align}

Collecting the estimates (\ref{j51}), (\ref{i1}) and (\ref{i21})--(\ref{j52}), choosing 
$\gamma>0$ small enough, and invoking the stability estimate (\ref{stability}), we can finally
conclude from Gronwall's lemma that 
\begin{equation}\label{esdiff}
\|\yl\|_{C^0([0,T]; L^2(\Omega))\cap L^2(0,T;
  H^1(\Omega))}^2\,\leq\,C_{13}\Big(\int_0^T
\|{\bf h}(t)\|^2_{L^3(\Omega)^3}
\dt\Big)^2\,\leq\,C_{14}\,\|{\bf h}\|_{\mathcal{V}}^4\,.
\end{equation}
This concludes the proof of Proposition~\ref{prop-differ}. \qed 

\vspace{6mm}
Using the convexity of ${\cal V}_{\rm ad}$, we immediately conclude from Proposition~\ref{prop-differ} the following 
result.

\vspace{3mm}
\beco\label{co-var}
Assume that the hypotheses {\bf (H1)}--{\bf (H8)} are fulfilled, and let $\bar \vb\in \mathcal{V}_{\rm ad}$ be an optimal control for problem {\bf (CP)} with associated state $\bar\vp={\cal S}(\bar \vb)$. Then we have for every ${\bf v}\in \mathcal{V}_{{\rm ad}}$ the inequality
\begin{align}\label{ineq-var}
&\beta_1\int_0^T\!\!\int_\Omega (\bar\varphi-\varphi_Q)\,\xi^{\bf h}\dx\ds \,+\,\beta_2\int_\Omega(\bar\vp(T)-\vp_\Omega)
\,\xi^{\bf h}(T)\dx\\
\no
&+\,\beta_3\int_0^T\!\!\int_\Omega \bar \vb\cdot({\bf v}-\bar \vb)\dx\ds\,\geq \,0\,,
\end{align}
where $\xi^{\bf h}$ is the unique solution to the linearized system {\rm (\ref{l1})--(\ref{l4})} associated with 
${\bf h}={\bf v}-\bar \vb$. 
\enco

\paragraph{The adjoint system and first-order necessary optimality conditions} In order to 
establish the necessary first-order optimality conditions for {\bf (CP)}, we need to eliminate $\xi^{\bf h}$ from inequality \eqref{ineq-var}. 
To this end, we introduce the {\em adjoint system} which formally reads as follows:
\begin{align}\label{ad1}
&-p_t-c_0\,\Delta p-\nabla p\cdot\Big[\bar \vb+m'(\bar\vp)\nabla \Big(\int_\Omega k(|x-y|)(1-2\bar\vp(y,t))\dy
\Big)\Big]\\[1mm]
\no
&-2\int_\Omega \nabla p(y,t)\,m(\bar\vp(y,t))\cdot\nabla k(|x-y|)\dy = \beta_1(\bar \vp-\vp_Q)\quad \hbox{ in }\,Q\,,\\
\label{ad2}
&\frac{\partial p}{\partial {\bf n}}=0\quad\hbox{on }\,\Sigma\,,\\
\label{ad3}
&p(T)=\beta_2(\bar\vp(T)-\vp_\Omega)\quad\hbox{a.e. in }\,\Omega\,.
\end{align}

\noindent Since the final value $p(T)$ only belongs to $L^2(\Omega)$, we can at best expect the regularity 
\[ p\in  H^1(0,T; H^1(\Omega)^*)\cap C^0([0,T]; L^2(\Omega))\cap L^2(0,T; H^1(\Omega))\,,\]
which entails that \eqref{ad1}--\eqref{ad2} must be understood in the weak variational sense. To this end,
we rewrite (\ref{ad1})--(\ref{ad2}) in the form 
\begin{align}\label{ad1bis}
&\duav{p_t(t), \eta}\,+\,c_0\int_\Omega\nabla p(t)\cdot\nabla \eta\dx \\[1mm]
\no
&-\int_\Omega \eta \,\nabla p(t) \cdot\Big[\bar\vb(t)\,+\,m'(\bar\vp(t))
\nabla \Big(\int_\Omega k(|x-y|)(1-2\bar\vp(y,t))\dy\Big)\Big]\dx\\[1mm]
\no
&-2\int_\Omega\! \eta \!\int_\Omega\!
 \nabla p(y,t)\,m(\bar\vp(y,t))\cdot\nabla k(|x-y|)\dy\dx\,=\,\int_\Omega\! \eta \,
 \beta_1\,(\bar \vp(t)-\vp_Q(t))\dx\,, 
\end{align}
for every $\eta\in H^1(\Omega)$ and almost every $t\in (0,T)$.   
 
We have the following existence and uniqueness result.
\bepr\label{exi-ad}
The adjoint system {\rm (\ref{ad1})--(\ref{ad3})}, written in the weak form {\rm (\ref{ad1bis})},  has a unique solution 
\[ p \in H^1(0,T; H^1(\Omega)^*)\cap C^0([0,T]; L^2(\Omega))\cap L^2(0,T; H^1(\Omega))\,.\]
\empr
\proof The proof is analogous to the first part of the proof of Proposition~\ref{prop-exilin}. 
In fact, one can again devise a Faedo-Galerkin approximation scheme, for which estimates
similar to the ones leading to (\ref{j41}) can be performed. An estimate resembling (\ref{j49})
cannot be derived since $\,\beta_1\,(\bar\vp(T)-\vp_\Omega)\,$ does not necessarily belong to
$H^1(\Omega)$.  One then obtains a weak solution 
that enjoys the asserted regularity and turns out to be unique. Since these arguments are
rather standard and straightforward, we can 
allow ourselves to omit the details here. \qed    

\vspace{5mm}
We are now in the position to eliminate $\,\xi^{\bf h}\,$ from (\ref{ineq-var}). We have the
following result. 

\bete\label{th-nec}
Assume that the hypotheses {\bf (H1)}--{\bf (H8)} are fulfilled, and let $\bar \vb\in \mathcal{V}_{\rm ad}$ be an optimal control for problem {\bf (CP)} with associated state $\bar\vp={\cal S}(\bar \vb)$ and adjoint state $p$. Then we have for every ${\bf v}\in \mathcal{V}_{{\rm ad}}$ the inequality
\begin{equation}\label{proj}
\beta_3\int_0^T\!\!\int_\Omega \bar \vb\cdot({\bf v}-\bar \vb)\dx\dt\,+\,\int_0^T\!\!\int_\Omega p({\bf v}-\bar \vb)\cdot\nabla\bar\vp\dx\dt\,\geq\, 0\,.
\end{equation}
\ente
\proof This is a standard calculation that can be left to the reader. We only note that we have
\begin{eqnarray}
\hspace*{1cm}&&\beta_1\int_0^T\!\!\int_\Omega (\bar\varphi-\varphi_Q)\,\xi^{\bf h}\dx{\dt} \,+\,\beta_2\int_\Omega(\bar\vp(T)-\vp_\Omega)
\,\xi^{\bf h}(T)\dx\\[2mm]
\hspace*{1cm}&&=\beta_1\int_0^T\!\!\int_\Omega (\bar\varphi-\varphi_Q)\,\xi^{\bf h}\dx{\dt} \,+\,
\int_0^T 
\left(\langle p_t(t),\xi^{\bf h}(t)\rangle\,+\,\langle\xi^{\bf h}_t(t),p(t)
\rangle\right) %\dx
\dt\nonumber\\[2mm]
\hspace*{1cm}&&=\int_0^T\!\!\int_\Omega p\,({\bf v}-\bar\vb)\cdot\nabla\bar\vp\dx\dt\,,\nonumber
\end{eqnarray}
where the last equality easily follows from expressing 
$\,p_t(t)\,$ and $\,\xi^{\bf h}_t(t)$  via the adjoint
equation (\ref{ad1bis}) and the linearized system (\ref{l1})--(\ref{l4}), and then integrating by parts. \qed

\vspace{7mm}
\noindent
{\bf Remark 5.}  \,\,The state system (\ref{p1})--(\ref{p5}), written for $\vp=\bar\vp$, the adjoint system
and the variational inequality (\ref{proj}) form together the first-order necessary optimality
conditions. 
Observe that we have $\,p\in L^2(0,T;H^1(\Omega))\cap
C^0([0,T];L^2(\Omega))\,$ and $\,\bar\vp\in L^\infty(0,T;H^2(\Omega))$, whence it follows that $\,p\,\nabla\bar\vp\in L^2(Q)^3\cap L^\infty(0,T;L^{3/2}(\Omega)^3)$, so that 
the variational inequality (\ref{proj}) is meaningful. Moreover, since $\mathcal{V}_{\rm ad}$ is a nonempty, closed, and convex subset of $L^2(Q)^3$, we
can infer from \eqref{proj}  
that for $\beta_3>0$ the optimal control $\bar \vb$ is the $L^2(Q)^3$-orthogonal projection of 
$-\beta_3^{-1}p\,\nabla\bar\vp$ onto $\mathcal{V}_{\rm ad}$. In particular, if the function $\,\widetilde{\bf v}=(\widetilde{v}_1,\widetilde{v}_2,
\widetilde{v}_3)\in L^2(Q)^3$, which is given by
\begin{equation}
\widetilde {v}_i(x,t):=\max\,\left\{\widetilde{v}_{1_i}(x,t),
\,\min\,\left\{\widetilde v_{2_i}(x,t),\,-\beta_3^{-1}\,p(x,t)\,
\partial_i\bar\varphi(x,t)\right\}\right\},
\end{equation}
for $i=1,2,3,$ and almost every $(x,t)\in Q$, belongs to 
${\cal V}_{\rm ad}$, then $\,\widetilde {\bf v}=\bar {\bf v}\,$, and 
the optimal control $\bar{\bf v}$ turns out to be a pointwise projection.
Notice, however, that the requirement $\widetilde{\bf v}\in {\cal V}_{\rm ad}$ implies that we should have $\widetilde{\bf v}_t\in L^2(0,T;L^3(\Omega)^3)$, which in general cannot be expected
since we only can guarantee the regularity $\,p_t\in L^2(0,T;H^1(\Omega)^*)$.
Therefore, the information about the optimal control that can be
recovered from the projection property may be rather weak, in general. This is in contrast to the 
non-convective local case (see, e.g., \cite[Thm. 3.16]{hw}) and to the
convective local 2D case (see \cite{zl14}, where different boundary conditions
are considered); it is in fact the price to be paid for considering the
three-dimensional case with the flow velocity as the control parameter.  

\vspace{5mm}
{\em Acknowledgement.} \,\,We thank an anonymous referee for his helpful 
comments on the first version of this paper.  
%%%%%%%%%%%%%%%%%%%%%%%%%%%%%%%%%%%%%%%%%%%
%%%%%%%%%%%%%%%%%%%%%%%%%% biblio %%%%%%%%%%%%%%%%%%%%%%%%%%%%%%%%%%%%%%%%%%%

\end{document}